\documentclass[a4paper]{amsart}
\usepackage{amsmath, amssymb, amsfonts}

\usepackage{enumerate}
\usepackage[all]{xy}

\newcommand{\cha}{\mathrm{char}\,}
\newcommand{\Spec}{\mathrm{Spec}\,}
\newcommand{\Sp}{\mathrm{Sp}\,}
\newcommand{\Spa}{\mathrm{Spa}\,}
\newcommand{\Cont}{\mathrm{Cont}\,}
\newcommand{\cont}{\mathrm{cont}\,}
\newcommand{\Spf}{\mathrm{Spf}\,}
\newcommand{\R}{\mathfrak{R}}
\newcommand{\QQ}{\mathbb{Q}}
\newcommand{\ZZ}{\mathbb{Z}}
\newcommand{\ad}{\mathrm{ad}}
\newcommand{\an}{\mathrm{an}}
\newcommand{\co}{\mathrm{co}}
\newcommand{\tr}{\triangleright}
\newcommand{\cat}[1]{\textrm{(#1)}}
\newcommand{\BERK}{\cat{$\Phi_K$-An}}
\newcommand{\KAN}{\cat{$K$-An}}

\DeclareMathOperator{\id}{id} 
\DeclareMathOperator{\Hom}{Hom}

\usepackage{amscd, amsthm}

\title{Weil restriction of $p$-adic analytic spaces}
\author{Christian Wahle}

\numberwithin{equation}{subsection}
\newtheorem{proposition}{Proposition}[subsection]
\newtheorem{corollary}[proposition]{Corollary}
\newtheorem{theorem}[proposition]{Theorem}
\newtheorem{lemma}[proposition]{Lemma}

\theoremstyle{definition}
\newtheorem{definition}[proposition]{Definition}
\newtheorem{example}[proposition]{Example}

\begin{document}
\bibliographystyle{amsplain}
	\maketitle
	\begin{abstract}
				We study the functor of Weil restriction in the category of Huber's adic spaces and in the category of Berkovich spaces.
		We prove criteria for the representability of these functors in the respective categories.
		As an application of adic Weil restrictions, we prove that adic N\'eron models behave well with respect to tamely ramified base change.
	\end{abstract}

	\section{Introduction}
	Weil restriction constitutes a fundamental tool in many different branches of geometry.
For instance, the Deligne torus $\mathbb S$, whose representations correspond to real hodge structures, is the Weil restriction of the complex torus to the reals.
Weil restriction also defines the notion of a jet space, a central object in motivic integration theory.
And, finally, Weil restriction is an important technique in arithmetic algebraic geometry, as it helps to understand the base change of N\'eron models.

Weil restriction is defined to be the adjoint of the base change functor.
The most important question related to the functor of Weil restriction is, of course, the question of representability of this functor.
In algebraic geometry, Weil restrictions of varieties were first studied by its inventor A.~Weil using Galois descent in the case where the base is a field (cf. \cite{weil}).
A detailed and more general treatment of the representability in the case of schemes can be found in \cite{blr}, 7.6.
The formal and rigid Weil restriction is discussed in \cite{bert}.

In this paper, we study the functor of Weil restriction for Berkovich spaces and Huber's adic spaces.
In both cases we first prove the existence of Weil restriction of affinoid spaces of finite type with respect to a finite and free morphism of base change of rank $n$ (cf. proposition \ref{restriktion_affinoider_adischer_raeume} resp. proposition \ref{berko_restriktion_affinoid}).
We obtain these results by constructing the Weil restrictions of discs of relative dimension $1$ (cf. proposition \ref{eindimensionaler_fall_darstellbar} resp. proposition \ref{berko_kreisscheibe_res}) and by observing that Weil restriction is compatible with products and closed immersions.
The Weil restriction of a disc is a (not necessarily quasi-compact) analytic domain in the analytification of the affine $n$-space.

In the non-affinoid case we first show that the restrictions of the affinoid subdomains can be glued together.
Our main results, theorem \ref{adisch_darstellbar} resp. theorems \ref{berko_darstellbar_gut} and \ref{berko_schlechte_res_darstellbar}, provide sufficient criteria for the result of this glueing process to represent the functor of Weil restriction.
Furthermore we can generalize the result of \cite{bert} on the representability of rigid Weil restrictions using adic spaces (cf. theorem \ref{theorem_pad}).
We show that Weil restriction is compatible with adification of schemes and rigid spaces (cf. proposition \ref{vergleich_adische_algebraische_weilres}).

As an application we use our results on the Weil restrictions of adic spaces in order to provide an adic analogue of the result of Edixhoven on the tamely ramified base change of N\'eron models (cf. \cite{be}, 4.2).
Namely, given a smooth and separated adic space $X_K$ over a discretely valued field $K$, corollary \ref{tame_base_change} states that, if there exists an adic N\'eron model $X'$ after a tamely ramified field extension, then an adic N\'eron model over the base ring can be constructed by taking the adic closure of $X_K$ in the Weil restriction of the N\'eron model $X'$.
The main argument in the proof of this result is the fact that this closure can be realized as the space of fixed points of a suitable Galois action on the Weil restriction of $X'$, and that this space of fixed points is smooth (cf. proposition \ref{fixpunkte_glatt}).

	\section{Weil restriction of adic spaces}	

	\subsection{Generalities about adic spaces}
The standard references to Huber's theory of adic spaces are \cite{hu1}, \cite{hu2} and \cite{hu3}.
An \emph{affinoid ring} is a pair $(A,A^+)$ consisting of an f-adic ring $A$ and a ring of integral elements $A^+$.
A morphism of affinoid rings is a homomorphism of rings respecting the corresponding rings of integral elements.
We will often write $A$ instead of $(A,A^+)$, whenever it is clear which ring of integral elements we want to use.
An \emph{f-adic ring} is a topological ring possessing an open subring $A_0$ which, in the induced topology, is adic\footnote{By this we mean that the ring topology coincides with the $\mathfrak a$-adic topology for some ideal of definition $\mathfrak a \subset A_0$; we do not require $A_0$ to be complete or separated.}, with a finitely generated ideal of definition; such a subring $A_0$ is called a \emph{ring of definition}.
When we refer to an affinoid ring $A$ we always assume that $A$ is a strictly noetherian Tate-ring, or that $A$ has a noetherian ring of definition; in both cases it is known that the presheaf of adic functions on $\Spa A$ is in fact a sheaf.

If $K$ is a complete non-archimedean field (which we will always assume to be non-trivially valued), we regard $K$ as an affinoid ring using the valuation ring $R \subset K$ as a ring of integral elements.
If $A$ is an adic ring (not necessarily complete and separated) with finitely generated ideal of definition, we regard $A$ as an affinoid ring, taking $A$ itself as a ring of integral elements.

Quotients of complete affinoid rings with respect to finitely generated ideals are again complete:
\begin{lemma}\label{quotientensindvollstaendig}
Let $A$ be a complete and separated f-adic ring and let $\mathfrak a \subset A$ be a finitely generated ideal.
Then $A/\mathfrak a$ is complete and separated.
\end{lemma}
\begin{proof}
We first observe that a sequence $(c_n)_{n \in \mathbb N}$ in an f-adic ring is Cauchy iff the sequence $(c_{n+1} - c_n)_{n \in \mathbb N}$ tends to zero.

Let $A_0 \subset A$ be a ring of definition and let $I$ be an ideal of definition of $A$.
Furthermore, let $(A/\mathfrak a)_0$ be a ring of definition of $A / \mathfrak a$.
We first assume that $A$ is a strictly noetherian Tate-ring.
By \cite{hu1}, 2.2.4, $A$ has a countable system of neighborhoods of zero.
Then \cite{hu1}, 2.2.14 implies that $A/\mathfrak a$ is separated.
If $A$ has a noetherian ring of definition, \cite{hu1}, 2.3.33 (iii) yields the separatedness of $A/\mathfrak a$ in this case.

Let $(c_n)_{n \in \mathbb N}$ be a Cauchy sequence $A/\mathfrak a$.
Then the sequence $(d_n)_{n \in \mathbb N}$, defined by
$$
d_0 := c_0, \quad d_n := c_n - c_{n-1} \text{ for }n>1,
$$
tends to zero.
Hence there is an integer $n_0 \ge 0$ such that $d_n \in (A/\mathfrak a)_0$ holds for all $n \ge n_0$.
Let $a_0, \ldots, a_{n_0} \in A$ be lifts of $c_0, \ldots,c_{n_0}$.
Then there exists an unbounded ascending sequence $\mu (\nu), \nu \ge n_0$ such that $d_\nu \in \pi (I^{\mu ( \nu )})$ for all $\nu \ge n_0$.
Assume $c_\nu$ has already been lifted to $a_\nu$ for $\nu \ge n_0$, and choose $a_\nu' \in \pi^{-1} (d_{\nu+1}) \cap I^{\mu(\nu+1)}$.
Then $a_{\nu+1} := a_\nu + a_\nu'$ is a lift of $c_{\nu+1}$ and $(a_\nu')_{\nu \in \mathbb N}$ tends to zero.
It follows that $(a_n)_{n \in \mathbb N}$ is Cauchy with limit $a$, and, by continuity, that $\pi(a)$ is the limit of $(c_n)_{n \in \mathbb N}$.
\end{proof}

If $A$ is an affinoid ring, a finite subset $M \subset A$ will be called a \emph{radius with values in $A$}, if $A \cdot M$ is an open subset of $A$.
Now, if $(M_1,\dots,M_d)$ is a finite system of radii $M_i$ with values in $A$ (also called a \emph{polyradius with values in $A$}), we can consider the subring $A \langle x_1,\ldots,x_d \rangle_{M_1,\ldots, M_d} \subset \widehat A [[x_1,\ldots,x_d]]$, consisting of all formal power-series $\sum_{\nu \in \mathbb N^d} a_\nu x_1^{\nu_1}\cdots x_d^{\nu_d}$ with coefficients in the f-adic completion $\widehat A$ of $A$, such that for all neighborhoods $U$ of $0 \in A$ there exists a $\nu' \in \mathbb N^d$ such that $a_{\nu} \in M_1^{\nu_1} \cdots M_d^{\nu_d} \cdot U$ for all $\nu \in \mathbb N^d$ satisfying $\vert \nu \vert \ge \vert \nu' \vert$.
This ring can be canonically endowed with the structure of a complete affinoid ring, and it is universal in the category of complete affinoid $A$-algebras with continuous morphisms with respect to the property that $m_i x_i$ is contained in the ring of integral elements $A\langle x_1,\ldots,x_d \rangle_{M_1,\ldots,M_d}^+$ (cf. \cite{hu3}, 3.5).
We will refer to affinoid $A$-algebras of this type as \emph{generalized Tate-algebras}.

For any scheme $X$ locally of finite type over some base scheme $\Spec A$, $A$ an affinoid ring, we denote by $X^\ad$ the associated adic space over $\Spa A$ (cf. \cite{hu3}, 3.8).
 
	\subsection{The affinoid case}

Let us at first recall the formal definition of direct images of functors.

\begin{definition}
Let $C$ be a category and let $h \colon S' \longrightarrow S$ be a morphism such that fibred products with $h$ always exist.
Then the \emph{direct image} of a contravariant functor $F$ on $C_{S'}$, the category of $S'$-objects, is defined as
$$
h_* F := F (  \cdot \times_S S').
$$
If $F = \Hom_{S'}( \cdot , X')$ for some $X'/S'$, we will write $\R_{S'/S}(X'):=h_*F$.
This functor is called the \emph{Weil restriction} of $X'$ with respect to $h$.
\end{definition}

If $\R_{S'/S}(X')$ is representable by an $S$-object $\R$, we denote by 
$$
\Psi \colon \R \times_S S' \longrightarrow X'
$$
the canonical morphism.
It is immediate to see that direct images commute with base change and products (as long as all occurring fibred products exist).

The category of Berkovich spaces admits fibred products.
In the adic setting, since we always assume the morphism of base change $h \colon S' \longrightarrow S$ to be finite, the base change of any adic space $X$ over $S$ with respect to $h$ also exists.

In this section, let $\varphi \colon A \longrightarrow B$ be a finite, free morphism of rank $n$ of complete affinoid rings.
We set $S':= \Spa B$ and $S:= \Spa A$.
Let $N$ be a radius with values in $B$.
We will show that the Weil restriction of the relative disc $\Spa B \langle x \rangle_N$ can be represented by an analytic subset of the analytic affine space $(\mathbb A^n_{\Spec A})^\ad$.
For the proof we need a generalization of lemma 1.6 in \cite{bert}.

\begin{definition}
Let $(L, \vert \cdot \vert)$ be a valued field with valuation ring $L^+$ and value group $\Gamma$.
We define the \emph{spectral value $\sigma_{L^+} (p)$} of a monic polynomial 
$$
p(z) = z^n +c_1 z^{n-1} + \ldots + c_n \in L[z]
$$
with respect to $L^+$ to be
$$
\sigma_{L^+} (p) := \max_{i=1,\ldots,n} \vert c_i \vert_{L^+}^{1/i} \in \Gamma \otimes_\ZZ \QQ \cup \{ 0\}.
$$
\end{definition}

As in \cite{bgr}, 1.5.4/1 one shows that the spectral value satisfies the equation
$$
\sigma_{L^+}(pq) = \max \{\sigma_{L^+}(p), \sigma_{L^+}(q)\}
$$
for all monic $p,q \in L[z]$.

\begin{lemma}\label{charpol_ganz_kriterium}
For any $b \in B$ let $\chi_b (z) \in A [z]$ denote the characteristic polynomial of $b$ over $A$.
Then we have $b \in B^+$ if and only if $\chi_b(z) \in A^+[z]$.
\end{lemma}
\begin{proof}
The ``if''-part is clear, since our assumption on $\varphi : A \longrightarrow B$ being finite implies that $B^+$ is the integral closure of $A^+$ in $B$.

Now assume $b \in B^+$, and furthermore assume that $A$ is an affinoid field.
Since $\varphi$ is finite, $b$ is the root of a monic polynomial $h \in A^+[z]$.
Let $g$ be the minimal polynomial of $b$ over $A$.
Then $1 \ge \sigma_{A^+} (h) \ge \sigma_{A^+} (g)$ implies $g \in A^+[z]$, hence $\chi_b(z) \in A^+[z]$.

Now let $A$ be an arbitrary affinoid ring, and let $x \in \Spa A$.
If $k(x)$ denotes the residue field of $\Spa A$ at $x$, we get a cocartesian diagram
$$
\xymatrix{
B \otimes_A \widehat{k(x)} & B \ar[l] \\
\widehat{k(x)}\ar[u]^{\varphi'}  &A \ar[l]^{\psi_x}\ar[u]^\varphi,
}
$$
of affinoid rings, where $\varphi'$ is finite and free, and where $\psi_x$ is the composition of the canonical projection with the completion morphism $\gamma \colon k(x) \longrightarrow \widehat{k(x)}$.

The characteristic polynomial $\chi_{b \otimes 1} (z) \in \widehat{k(x)}[z]$ equals the image of $\chi_b(z)$ under $A[z] \longrightarrow \widehat{k(x)}[z]$.
Since
$$
b \otimes 1 \in B^+ \otimes_{A^+} \widehat{k(x)^+} \subset (B \otimes_A \widehat{k(x)})^+,
$$
we can conclude  $\chi_{b \otimes 1} \in \widehat{k(x)^+} [z]$.
Thus, it suffices to verify the following assertion:\\

\emph{
$(\ast)$\quad For $f \in A$ and $x \in \Spa A$ write $\bar f_x := \psi_x (f)$.
If $\bar f_x \in \widehat {k(x)}^+$ holds for all $x \in \Spa A$, then $f \in A^+$.
}
\\

Let $f \in A$ be such an element.
It suffices to show that $f \in A^+$, given $\vert f \vert_x \le 1$ for all $x \in \Spa A$.
Using \cite{hu2}, 3.3 it is even sufficient to check that, whenever $x \in \Cont A$ fulfills $\vert b \vert_x \le 1$ for all $b \in A$, then $\vert b \vert_x \le 1$ is true for all elements $b$ in the integral closure of $A^+[f]$ in $A$.

But if $b \in A$ is the root of a polynomial
$$
z^n + p_1 z^{n-1} + \ldots + p_n \in (A^+[f])[z],
$$
then we must have $\vert p_i \vert_x \le 1$ for all $i = 1,\ldots,n$, and, hence, $\vert b \vert \le 1$.
\end{proof}

We now construct the adic space over $\Spa A$ representing $\R_{S'/S}(\Spa B \langle x \rangle_N)$.
Let us at first choose a free system $e_1,\ldots,e_n \in B$ of generators over $A$ and write $N =\{r_1,\ldots,r_k\}$.
Let $M\subset A^{\circ \circ}$ be a finite and non-empty set of topologically nilpotent elements.
For all integers $\lambda \ge 0$ we denote by $M(\lambda)$ the set of $\lambda$-fold products of elements in $M$.
We define
$$
A_\lambda := A \langle X,Y_1,\ldots,Y_k \rangle_{\underbrace{M(\lambda), \ldots ,M( \lambda )}_n,\underbrace{\{1\},\ldots, \{1\}}_{nk}}
$$
for systems $X=(x_1,\ldots,x_n)$ and $Y_i = (y_{i1},\ldots,y_{in})$, $i=1,\ldots,k$ of indeterminates.
For all $1 \le i \le k$ let $c_1(r_iX), \ldots, c_n (r_iX) \in A[X]$ be the coefficients of the characteristic polynomial of $r_i \sum_{j=1}^n x_j e_j \in B [X]$ and let $I_N \subset A [X,Y_1,\ldots,Y_k]$ be the ideal generated by
$$
y_{ij} - c_j(r_iX), \qquad i=1,\ldots,k, \quad j=1,\ldots,n.
$$
For all $\lambda \in \mathbb N$ we put 
\begin{equation}\label{c_lambda_definition}
C_{\lambda,N} = A_{\lambda} / I_N \cdot A_\lambda.
\end{equation}

\begin{lemma}\label{vertraeglich_offen}
The following assertions hold:
\begin{enumerate}[(i)]
\item For all integers $\lambda' \ge \lambda \ge 0$ there is a continuous morphism
$$
\rho^*_{\lambda,\lambda'} \colon C_{\lambda',N} \longrightarrow C_{\lambda,N}
$$
of complete affinoid rings over $A$, such that $\rho^*_{\lambda,\lambda''} = \rho^*_{\lambda,\lambda'} \circ \rho^*_{\lambda',\lambda''}$ for all integers $\lambda '' \ge \lambda ' \ge \lambda \ge 0$.
\item The induced morphism
$$
\rho_{\lambda,\lambda'} \colon \Spa C_{\lambda,N} \longrightarrow \Spa C_{\lambda',N}
$$
is an open immersion satisfying $\rho_{\lambda,\lambda''} = \rho_{\lambda',\lambda''} \circ \rho_{\lambda,\lambda'}$ for all integers $\lambda'' \ge \lambda' \ge \lambda \ge 0$.
\end{enumerate}
\end{lemma}
\begin{proof}
The assertion follows right from the universal property of generalized Tate algebras and from the fact that we can identify $\Spa C_{\lambda,N}$ via $\rho_{\lambda,\lambda'}$ with the rational domain 
$$
U_{\lambda'} = \{ y \in \Spa C_{\lambda',N} ; \vert m ( \pi_{\lambda'} (x_i)) \vert_y \le 1 \text{ for all } m \in M(\lambda), 1 \le i \le n\}
$$
in $\Spa C_{\lambda',N}$.
\end{proof}

Now we can prove our first representability result.

\begin{proposition}\label{eindimensionaler_fall_darstellbar}
The Weil restriction $\R_{\Spa B/\Spa A} (\mathbb{B}^1_{B,N})$ is representable by the adic space
$$
\R_{\Spa B / \Spa A} (\mathbb{B}^1_{B,N} ) = \varinjlim_{\lambda \in \mathbb N} \Spa C_{\lambda,N},
$$
obtained by glueing the $\Spa C_{\lambda,N}, \lambda \in \mathbb N$.
\end{proposition}
\begin{proof}
Let $\lambda \in \mathbb N$.
We define
$$
\Psi_\lambda^* \colon B \langle x \rangle_N \longrightarrow C_{\lambda,N} \otimes_A B
$$
by setting $\Psi_\lambda^*(x) = \sum_{j=1}^n x_j e_j$.
This yields a well-defined continuous morphism of complete affinoid rings, if
$$
r_i \cdot \sum_{j=1}^n x_je_j \in (C_{\lambda,N} \otimes_A B)^+
$$
holds for all $i=1,\ldots,k$.
But the latter is true, namely by construction of the $C_{\lambda,N}$ and by lemma \ref{charpol_ganz_kriterium}.

Now we can glue the $\Spa B$-morphisms $\Psi_\lambda := \Spa (\Psi_\lambda^*)$, $\lambda \in \mathbb N$ to obtain an $\Spa B$-morphism 
$$
\Psi \colon \left( \varinjlim_{\lambda \in \mathbb N} \Spa C_{\lambda,N} \right) \times_{\Spa A} \Spa B \longrightarrow \mathbb{B}^1_{B,N}.
$$
We show that this morphism defines a bijection
\begin{eqnarray*}
\Theta \colon \Hom_{\Spa A} (Z, \varinjlim_{\lambda \in \mathbb N} \Spa C_{\lambda,N}) & \longrightarrow & \Hom_{\Spa B} ( Z \times_{\Spa A} \Spa B, \mathbb B^1_{B,N}),\\
a&\longmapsto&\Psi \circ ( a \times \id),
\end{eqnarray*}
for any adic space $Z$ over $\Spa A$.
We may assume $Z =\Spa D$, where $D$ a complete affinoid ring.
Furthermore, since $\varphi : A \longrightarrow B$ is finite, the fibered product $Z \times_S S'$ exists and is affinoid, with ring of global sections $D \otimes_A B$.

Now consider $a_1,a_2 \in \Hom_{\Spa A} (Z, \varinjlim_{\lambda \in \mathbb N} C_{\lambda,N})$ such that $\Theta(a_1)=\Theta(a_2)$.
As $Z$ is quasi-compact, there is some $\lambda \in \mathbb N$ such that $a_1$ und $a_2$ factor through $\Spa C_{\lambda,N}$ and $\Psi_{\lambda} \circ (a_1 \times \id) = \Psi_{\lambda} \circ (a_2 \times \id)$.
Let $a_1^*$ (resp. $a_2^*$) denote the continuous morphisms of affinoid rings associated to $a_1$ (resp. $a_2$).
Then
\begin{eqnarray*}
\sum_{j=1}^n a_1^*(x_j) e_j &=& (a_1^* \otimes \id) (\sum_{j=1}^n x_j e_j)\\
&=&( \Psi_\lambda \circ (a_1 \times \id))^*(x)\\
&=&( \Psi_\lambda \circ (a_2 \times \id))^*(x)\\
&=&\sum_{j=1}^n a_2^*(x_j) e_j,
\end{eqnarray*}
hence $a_1 = a_2$.

Now consider an $\Spa B$-morphism $\tau \colon Z \times_{\Spa A} \Spa B \longrightarrow \mathbb B^1_{B,N}$ and denote the associated continuous morphism of complete affinoid rings by $\tau^*$.
We write 
$$
\tau^*(x) = \sum_{j=1}^n a_j e_j,\qquad a_1,\ldots,a_n \in D.
$$
Now, since $M$ consists of a finite number of topologically nilpotents, there exists a $\lambda_0 \in \mathbb N$ such that $m \cdot a_j \in D^+$ for all $j=1,\ldots,n$ and all $m \in M(\lambda), \lambda \ge \lambda_0$.
From
$$
r_i \sum_{j=1}^n a_j e_j \in (D \otimes_A B)^+ \quad \text{for all } i=1,\ldots,k
$$
we can conclude, using lemma \ref{charpol_ganz_kriterium}, that all coefficients $c_1 (r_i \tau^*(x)), \ldots, c_n (r_i \tau^*(x)) \in D$ of the characteristic polynomial $r_i \sum_{j=1}^n a_j e_j \in D \otimes_A B$ are in fact contained in $D^+$ for all $i =1,\ldots,k$.
Hence, for all $\lambda \ge \lambda_0$ there exists an $A$-morphism
$$
\psi_\lambda \colon A_\lambda \longrightarrow D
$$
satisfying $\psi_\lambda(x_j) = a_j$ and $\psi_\lambda(y_{ij}) = c_j(r_i \tau^*(x))$ for all $i=1,\ldots,k$, $j=1,\ldots,n$.
This morphism factors as
$$
\xymatrix{
A_\lambda \ar[r]^{\psi_\lambda} \ar[d]&D\\
C_{\lambda,N} \ar[ru]_{a^*}&
}
$$
and finally $\tau^* = (a^* \otimes \id) \circ \Psi_\lambda^*$.
\end{proof}

The construction suggests that, in general, the Weil restriction of a unit disc might in general not be quasi-compact.
Indeed, the following example shows that such a situation can already occur in the rigid, and therefore also in the adic setting, if we study the Weil restriction with respect to a finite field extension $K'/K$.
Note that, due to the fact that Weil restriction is compatible with products, it follows immediately that $\R_{K'/K}(\mathbb B^1_K)$ is quasi-compact given $K'/K$ is separable.

\begin{example}\label{non-quasi-compact-example}
Put $K:=\mathbb F_p (x)$, $K' := K (\sqrt[p]{x}) = K[t]/(t^p-x)$.
We consider the non-archimedean valuation $\vert \cdot \vert$ on $K$, which satisfies $\vert f/g \vert = r^{\deg f - \deg g}$, where $0<r<1$.
Then $K' \otimes_K K' = K'[y] / (y-\bar t)^p$.
For any integer $k \ge 0$, let $b_k$ be the image of $x^{-k} (y- \bar t)$ in $K' \otimes_K K'$.
Then all $b_k$ are nilpotent, i.e. $\vert b_k \vert_{\sup} =0$.
However, $\vert x^{-k} \vert_{\sup} \to \infty$, hence the Weil restriction of $\mathbb B^1_{K'}$ with respect to $K'/K$ cannot be quasi-compact.
\end{example}

The compatibility of the Weil restriction functor with products implies the representability of the Weil restriction of higher dimensional discs:

\begin{proposition}\label{polykreis_darstellbar}
Let $N =(N_1 , \ldots, N_d)$ be a polyradius with values in $B$, $d \ge 1$.
Then the Weil restriction of 
$$
\mathbb B^d_{B,N} = \Spa B \langle x_1,\ldots,x_d \rangle_{N_1,\ldots,N_d}
$$
exists and
\begin{eqnarray*}
\R_{\Spa B / \Spa A} (\mathbb B^d_{B,N}) &=& \prod_{i=1}^d \R_{\Spa B / \Spa A} ( \mathbb B^1_{B,N_i})\\
&=&\varinjlim_{\lambda \in \mathbb N} \Spa \widehat \bigotimes_{i=1,\ldots,d} C_{\lambda,N_i}.
\end{eqnarray*}
\end{proposition}
\label{ad_res_affinoide}

We now prove the representability of $\R_{S'/S} (X')$ where $X'$ is an affinoid adic space of finite type over $S'$, which means that there exists a quotient morphism
\begin{equation}
\pi \colon B \langle X \rangle_{N} \longrightarrow D
\end{equation}
of complete affinoid rings, where $X' = \Spa D$, $X =(x_1,\ldots,x_d )$ and where $N=(N_1,\ldots,N_d)$ is a polyradius with values in $B$.

\begin{proposition}\label{restriktion_affinoider_adischer_raeume}
The Weil restriction $\R_{S' / S} (X')$ of an affinoid adic space $X' = \Spa D$ is an adic space, locally of finite type over $S$.
\end{proposition}
\begin{proof}
We may assume that $D$ is complete and we fix an $A$-basis $e_1,\dots,e_n$ of $B$.
Let $\mathfrak a$ be the kernel of the above quotient morphism $\pi$.
We put $E_{\lambda,N} := \widehat{\bigotimes}_{k=1}^d C_{\lambda,N_k}$ for every $\lambda \in \mathbb N$. (For the definition of the $C_{\lambda,N}$ cf.  \eqref{c_lambda_definition}.)
We denote the first $n$ coordinate functions of $C_{\lambda,N_k}$ by $\overline{x_{k1}}, \ldots, \overline{x_{kn}}$. 
For all $k=1,\ldots,d$ and all $\lambda \in \mathbb N$ we have the canonical morphism $\sigma_k \colon C_{\lambda,N_k} \longrightarrow E_{\lambda,N}$.
From the representability of the Weil restriction of polydiscs we obtain a canonical continuous $B$-morphism
$$
\Psi_\lambda^* \colon B \langle X \rangle_{N} \longrightarrow E_{\lambda,N} \otimes_A B, \qquad \lambda \in \mathbb N
$$
of complete affinoid rings, given by
$$
\Psi_\lambda^* (x_k) = \sum_{j=1}^n \sigma_k(\overline{x_{kj}}) e_j, \qquad k=1,\ldots,d.
$$
Now let $(f_i)_{i \in I}$ be a generating system of $\mathfrak a$.
For all $i \in I$ and all $\lambda \ge 0$ we can write
$$
\Psi_\lambda^* (f_i) = \sum_{j=1}^n f_{ij\lambda} e_j, \qquad f_{ij\lambda} \in E_{\lambda,N}.
$$
Let $\mathfrak a^\co_\lambda$ be the ideal in $E_{\lambda,N}$ generated by the functions $f_{i1\lambda}, \ldots, f_{in\lambda}$, $i \in I$. (We will refer to this as the \emph{ideal of coefficients of $\mathfrak a_\lambda$ in $E_{\lambda,N}$\label{koeffizientenideal}}.)

We now prove that the ideal $\mathfrak a^\co_\lambda$ defines the adic space representing the Weil restriction of $X'$.
Let $Z$ be an affinoid adic space over $\Spa A$ and let $\lambda \in \mathbb N$.
Consider an $A$-morphism $a^* \colon E_{\lambda,N} \longrightarrow \mathcal O_Z (Z)$.
Since $\mathfrak a^\co_\lambda$ is contained in the kernel of $a^*$ if and only if $\mathfrak a_\lambda$ is contained in the kernel of $(a^* \otimes \id) \circ \Psi_\lambda^*$, we obtain a commutative diagram
$$
\xymatrix{
\Hom_A^\cont (\widehat{E_{\lambda,N}/\mathfrak a^\co_\lambda}, \mathcal O_Z(Z)) \ar[r]^{\rho_\lambda} \ar[d] & \Hom_B^\cont (D, \mathcal O_Z(Z) \otimes_A B) \ar[d]\\
\Hom_A^\cont (E_{\lambda,N},\mathcal O_Z(Z)) \ar[r] & \Hom_B^\cont (B\langle X \rangle_N, \mathcal O_Z(Z) \otimes_A B).
}
$$
Now let $a \colon Z \longrightarrow \varinjlim_{\lambda \in \mathbb N} \Spa (E_{\lambda,N} / \mathfrak a^\co_\lambda)$ be an $\Spa A$-morphism.
Then there exists a $\lambda \in \mathbb N$ such that $a$ factors through $\Spa E_{\lambda,N} / \mathfrak a^\co_\lambda$.
Let $a^* \colon \widehat{E_{\lambda,N} / \mathfrak a^\co_\lambda} \longrightarrow \mathcal O_Z(Z)$ denote the associated continuous $A$-morphism of complete affinoid rings.
Then $\rho_\lambda (a^*)$ yields an $\Spa B$-morphism $Z \times_{\Spa A} \Spa B \longrightarrow \Spa D$, hence we have a map
$$
\Theta \colon \Hom_{\Spa A} (Z, \varinjlim_{\lambda \in \mathbb N} \Spa (E_{\lambda,N} /\mathfrak a^\co_\lambda)) \rightarrow \Hom_{\Spa B} (Z \times_{\Spa A} \Spa B, \Spa D).
$$
But now, using the quasi-compactness of affinoid adic spaces, the representability of the Weil restriction of polydiscs and the commuting diagram above, one can easily check that $\Theta$ is bijective.
\end{proof}

For $d=0$ the previous proposition shows that the restriction of a closed subspace $V \subset S'$ is a closed subspace in $S$.
But the same holds for open subspaces $V \subset S'$.
Namely, since finite morphisms of adic spaces are closed, we have the equation
$$
\R_{S'/S}(V) = S \setminus h(S' \setminus V).
$$
Hence, by the same functorial argument as in \cite{blr}, 7.6/2 we find that the formation of direct images of presheaves respects closed and open immersions.

We summarize what we have explained in the following lemma.
\begin{lemma}\label{adische_restriktion_vertraeglich_mit_offenen_immersionen}
Let $X'$ be an adic space over $S'$ such that the Weil restriction $\R_{S'/S}(X')$ exists as an adic space over $S$ and let $U' \hookrightarrow X'$ be an open immersion.
Then $\R_{S'/S}(U')$ can be represented as an open subspace of $\R_{S'/S}(X')$.
The canonical morphism $\Psi_{U'}$ associated to the Weil restriction of $U'$ is the restriction of the canonical morphism $\Psi_{X'}$ associated to the Weil restriction of $X'$.
\end{lemma}
\begin{proof}
We may write
$$
(\ast)\quad\R_{S'/S}(U') = \R_{S'/S}(U') \times_{\R_{S'/S}(X')} T = \R_{T'/T}(U' \times_{X'} T'),
$$
with $T:= \R_{S'/S}(X')$ and $T' := T \times_S S'$.
The morphism
$$
\R_{S'/S}(U') \times_S S' \longrightarrow  \R_{S'/S}(X') \times_S S' \longrightarrow X'
$$
allows us to view $\R_{S'/S}(U')$ as an $X'$-object.
From $(\ast)$ it follows that 
$$\
Psi_{U'} \colon \R_{S'/S}(U') \times_S S' = \R_{S'/S}(U') \times_T T' \longrightarrow U'
$$
is even an $X'$-morphism.
This in turn means $\Psi_{U'} = \Psi_{X'} \vert_{\R_{S'/S}(U') \times_S S'}$.
\end{proof}
	\subsection{A criterion for the global case}

In the following let $X'$ be an adic space locally of finite type over $S'$.
We consider the system $(U_i')_{i \in I}$ of all open affinoid subspaces of $X'$.
Since Weil restriction respects open immersions, we can now construct an adic space $\R$, locally of finite type over $S$, together with an $S'$-morphism $\R \times_S S' \longrightarrow X'$, by glueing the Weil restrictions $\R_{S'/S}(U_i')$, $i \in I$ with respect to the canonical glueing data.
However, similar to the case of (formal) schemes or rigid spaces, this space $\R$ will not always represent the Weil restriction of $X'$.

In this section we provide a sufficient condition for $\R$ to represent the Weil restriction of $X'$.

\begin{definition}
\begin{enumerate}[(i)]
\item An adic space $X$ has property (P) if any finite subset $I \subset X$ is contained in an open affinoid $U \subset X$.
\item Let $K$ be a complete non-archimedean field. A rigid $K$-space $X$ has property (PAd) if for any finite subset $I$ in the associated adic space $X^\ad$ there exists an open affinoid subspace $U \subset X$ such that $I \subset U^\ad$.
\end{enumerate}
\end{definition}

We see that if $X$ is a rigid $K$-space, which does not have property (PAd), its associated adic space $X^\ad$ might still satisfy property (P).

Before we give the proof of our main theorem on the Weil restriction of adic spaces, we need a technical lemma.

\begin{lemma}\label{technische_faktorisierung}
Let $X'$ be an adic space over $S'$ and let $U' \subset X'$ an open subspace.
If $\R_{S'/S}(X')$ is an adic space over $S$ and $a \colon Z \longrightarrow \R_{S'/S}(X')$ an $S$-morphism, such that the $S'$-morphism $a' \colon Z \times_S S' \longrightarrow X'$ corresponding to $a$ factors through $U'$, then $a$ factors through $\R_{S'/S}(U')$.
\end{lemma}
\begin{proof}
Assume $a'$ factors through $U'$.
Then there is an $S$-morphism 
$$
b \colon Z \longrightarrow \R_{S'/S}(U')
$$
with $\Psi_{U'} \circ (b \times \id_{S'}) = a'$, where we have seen in lemma \ref{adische_restriktion_vertraeglich_mit_offenen_immersionen} that the canonical morphism $\Psi_{U'}$ is the restriction of $\Psi_{X'}$ to $\R_{S'/S}(U') \times_S S'$.
Thus we have 
$$
\Psi_{X'} \circ (b \times \id_{S'}) = a' = \Psi_{X'} \circ (a \times \id_{S'}).
$$
Now the representability of $\R_{S'/S}(X')$ implies $a = b$.
In particular, the image of $a$ is contained in $\R_{S'/S}(U')$.
\end{proof}

\begin{theorem}\label{adisch_darstellbar}
Let $h \colon S' \longrightarrow S$ be a finite free morphism of affinoid adic spaces and let $X'$ be an adic space locally of finite type over $S'$ satisfying property (P).
Then the Weil restriction $\R_{S'/S} (X')$ exists as an adic space over $S$ and $\R_{S'/S}(X')$ is locally of finite type over $S$.
\end{theorem}
\begin{proof}
We have to prove that the adic space $\R$ we constructed above represents the Weil restriction $\R_{S'/S} (X')$.
Let $Z$ be an affinoid adic space over $S$ and put $Z'=Z \times_S S'$.
We consider $S$-morphisms $a_1,a_2 \colon Z \longrightarrow \R$ satisfying the equation 
$$
\Psi \circ (a_1 \times \id) = \Psi \circ (a_2 \times \id).
$$
Now let $x \in Z$.
We set $\tilde I := p^{-1} (x)$, where $p$ denotes the canonical projection $p \colon Z' \longrightarrow Z$.
By \cite{hu1}, 3.12.20 the set $\tilde I$ is finite.
Hence, by assumption on $X'$, there exists an open affinoid neighborhood $U' \subset X'$ of the set $I := \Psi \circ (a_1 \times \id) (\tilde I)= \Psi \circ (a_2 \times \id)(\tilde I)$, and $U'' := (\Psi \circ (a_1 \times \id))^{-1}(U') \subset Z'$ is an open neighborhood of $\tilde I$.
It follows that
$$
x \in \R_{Z \times_S S' /Z} (U'') = Z \setminus p (Z' \setminus U'') =:V_x
$$
and $V_x$ is open in $Z$.
The restrictions $a_1' := \Psi \circ (a_1 \times \id) \vert_{V_x \times_S S'}$ (and, hence, $a_2' := \Psi \circ (a_2 \times \id)\vert_{V_x \times_S S'}$) factor through $U' \subset X'$.
Once we can show that now both $a_1 \vert_{V_x}$ and $a_2 \vert_{V_x}$ factor through $\R_{S'/S}(U')$ we can conclude from the representability of $\R_{S'/S}(U')$ that $a_1 \vert_{V_x} = a_2 \vert_{V_x}$ holds.
Let $v \in V_x$.
We consider the morphism
$$
a_{1,v} \colon \Spa k(v) \longrightarrow V_x \longrightarrow \R.
$$
Let $W' \subset X'$ be an open affinoid subspace, which has the property that $\R_{S'/S}(W') \subset \R$ is an open neighborhood of the image of $v \in \Spa k(v)$ under $a_{1,v}$.
The points of $\Spa k(v)$ are generalizations of $v$.
It follows, that $a_{1,v}$ factors through $\R_{S'/S}(W')$.
But then, by definition of $\Psi$, we have 
$$
\Psi \circ ( a_{1,v} \times \id_{S'}) = \Psi_{W'} \circ (a_{1,v} \times \id_{S'}).
$$
Since the projections
$$
p_v \colon \Spa k(v) \times_S S' \longrightarrow \Spa k(v)
$$
are finite, we can use \cite{hu1}, 3.12.20 and \cite{hu4}, 1.3.4 in order to conclude that $p_v$ is specializing.
In particular, every point of $\Spa k(v) \times_S S'$ is a generalization of a point in $p^{-1}(v)$.
Now $v \in V_x$ implies that the image of $\Psi_{W'} \circ ( a_{1,v} \times \id_{S'})$ is contained in $U' \cap W'$.
Since $W'$ is affinoid, lemma \ref{technische_faktorisierung} shows that the image of $a_{1,v}$ is already contained in $\R_{S'/S}(U' \cap W') \subset \R_{S'/S}(U')$.
Thus $a_1 = a_2$.

Now let $a' \colon Z' \longrightarrow X'$ be an $S'$-morphism and let $x \in Z$. 
We put $I:= a'(p^{-1}(x))$.
Choose an open affinoid neighborhood $U_x' \subset X'$ of $I$.
Then $V_x \times_S S' \subset (a')^{-1}(U_x')$, where $V_x := \R_{Z'  / Z} ( (a')^{-1}(U_x'))$.
Thus $a'$ restricts to an $S'$-morphism $V_x \times_S S' \longrightarrow U_x'$, which in turn corresponds in a unique way to an $S$-morphism $a_x \colon V_x \longrightarrow \R_{S'/S}(U_x') \subset \R$.
Now it is straighforward to check that the $a_x$, $x$ varying over $Z$, glue to an $S$-morphism $a\colon Z \longrightarrow \R$ satisfying $\Psi \circ (a \times \id) = a'$.
\end{proof}

For later use we want to prove that property (P) is stable under taking Weil restriction.

\begin{proposition}\label{eigenschaft_p_steigt_ab}
Let $X'$ be an adic space, locally of finite type over $S'$.
If $X'$ satisfies property (P), then $\R_{S'/S}(X')$ satisfies property (P).
\end{proposition}
\begin{proof}
From the local construction we see that it suffices to reduce to the case where $X'$ is affinoid.
Let $I \subset \R_{S'/S}(X')$ be a finite set of points and let $\tilde I = p^{-1}(I)$, where now $p \colon \R_{S'/S}(X') \times_S S' \longrightarrow \R_{S'/S}(X')$ is the projection onto the first factor.
Since $\tilde I$ is finite, we find an open affinoid subspace $U' \subset X'$ containing $\Psi_{X'} (\tilde I ) \subset X'$.
Then we have $\tilde I \subset \Psi_{X'}^{-1} (U')$, and it suffices to show that $\tilde I$ is already contained in $\R_{S'/S}(U') \times_S S'$.
But again, using lemma \ref{technische_faktorisierung} and the fact that $p$ is specializing (cf. the argument in the previous theorem), it follows that $p^{-1}(x)$ is already contained in $\R_{S'/S}(U') \times_S S'$ for all $x \in I$, which settles the proof.
\end{proof}
	\subsection{Relation to algebraic and rigid geometry}
In this subsection we discuss the compatibility of Weil restriction with respect to analytification, completion and adification.

At first we should point out that the functor $t$, which associates to every locally noetherian formal scheme an adic space, does not commute with Weil restriction.
The reason is that the Weil restriction of an affinoid adic space of finite type does not need to be quasi-compact, whereas the Weil restriction of an affine locally noetherian formal scheme of topologically finite type is always quasi-compact.

Now let $K$ be a complete, non-archimedean valued field, let $R:=\{ x \in K ; \vert x \vert \le 1\}$ be the valuation ring of $K$ with residue field $k$.
We want to prove that property (PAd) is a sufficient condition for the representability of Weil restriction of rigid-analytic $K$-spaces.
First we check that property (PAd) is local on the base.

\begin{lemma}\label{PAdlokal}
Let $Y_K$ be a quasi-affinoid rigid $K$-space.
Then $Y_K$ has property (PAd).
\end{lemma}
\begin{proof}
The rigid space $Y_K$ is by definition quasi-compact and there exists an admissible formal $R$-scheme $X = \Spf A$ with generic fibre $X_K = \Sp A_K$, such that $Y_K \subset X_K$ is an admissible open.
Then there exists an admissible blow-up $\pi_\mathcal A \colon X_\mathcal A \longrightarrow X$, such that $X_\mathcal A$ contains an open formal subscheme $Y$, which is an admissible $R$-model of $Y_K$.
Let $I \subset Y_K^\ad$ be a finite set of points.
Then $Y_K^\ad \subset t(Y)$ and we denote the image of $I$ under the canonical morphism $\pi_Y \colon t(Y) \longrightarrow Y$ of topologically and locally ringed spaces by $I'$.
Since for any open formal subscheme $U \subset Y$ the equation $\pi_Y^{-1}(U) = t(U)$ holds, it suffices to see that there exists an open affine formal subscheme $U \subset Y$ such that $I' \subset U$.
But the latter is clear, since the special fibre $Y_k$ of $Y$ is quasi-projective over $A_k$.
\end{proof}

We can now use lemma \ref{PAdlokal}  in order to show that (PAd) is local on $S$.
Let $X$ be a rigid $K$-space over a rigid $K$-space $S$ and let $U \subset S$ be an open affinoid subspace.
Assume $X$ has property (PAd) and put $X_{U}:= X \times_{S} U$.
Let $I \subset (X_U)^\ad $ be a finite set of points.
Then there exists an open affinoid subspace $V \subset X$ such that $I \subset V^\ad$.
By shrinking $W$ we may assume that $W:= V \cap X_U$ is quasi-affinoid.
Lemma \ref{PAdlokal} yields the existence of an open affinoid subspace $W' \subset W$, such that $I \subset (W')^\ad$.

\begin{theorem}\label{theorem_pad}
Let $h \colon S' \longrightarrow S$ be a finite and locally free morphism of rigid $K$-spaces and let $X'$ be a rigid $K$-space over $S'$.
If $X'$ satisfies property (PAd), the Weil restriction $\R_{S'/S} (X')$ is a rigid $K$-space over $S$.
\end{theorem}
\begin{proof}
Let $\mathfrak{S}$ be an admissible open covering of $S$, such that for all $V \in \mathfrak{S}$ the preimage $h^{-1}(V)$ is $K$-affinoid and such that the induced homomorphism of $K$-affinoid algebras $\mathcal{O}_S (V) \longrightarrow \mathcal{O}_{S'}(h^{-1}(V))$ is finite and free.
Let $T$ be an affinoid $K$-space over $S$, $T':= T \times_S S'$, $\psi \colon T' \longrightarrow X'$ an $S'$-morphism and $\mathfrak{U}$ the covering of $X'$ consisting of all affinoid subspaces of $X' \times_S V$, where $V$ varies over $\mathfrak{S}$.
We want to use \cite{bert}, 1.14.
Thus we have to show that $(\R_{T'/T} (\psi^{-1}(U)))_{U \in \mathfrak{U}}$ is an admissible covering of $T$, because then the pull-back of this covering to $T'$ is a refinement of $\psi^* \mathfrak{U}$.

As we have seen, we may assume that both $S'$ and $S$ are affinoid, such that $\mathcal{O}_{S'}(S')$ is finite free over $\mathcal{O}_S(S)$.
The sets $\R_{T'/T} (\psi^{-1} (U)) \subset T$, $U \in \mathfrak{U}$ are admissible opens.
By \cite{hu1}, 3.4.17 we can check the admissiblity of a covering of $T$ by admissible opens by verifying that the associated adic open subspaces cover the adification $T^\ad$ of $T$.

Let $x \in T^\ad$.
The morphism $p^\ad \colon (T')^\ad \longrightarrow T^\ad$ is finite, hence  $\psi^\ad ( \tilde I ) \subset (X')^\ad$ is finite, where we have set $\tilde I := (p^\ad)^{-1}(x) \subset (T')^\ad$.
Now we choose a geometric open affinoid $U \subset (X')^\ad$, i.\,e. there exists and open affinoid subspace $U' \subset X'$ such that $U = (U')^\ad$ and  $\psi^\ad (\tilde I) \subset U$.
We put $V:= \R_{T'/T} (\psi^{-1}(U')) = T \setminus p(T' \setminus \psi^{-1}(U'))$, where $V$ is admissible open in $T$.

Consider the commutative diagram of ringed sites
\[\xymatrix{
\bar U'\ar@{^{(}->}[r]  & T'\ar[r]^p  &   T \\
\bar U\ar@{^{(}->}[r]\ar[u]&   T'^\ad \ar[r]^{p^\ad}\ar[u]^{\rho'}    &   T^\ad\ar[u]^\rho,
}\]
where we have set $\bar U':=\psi^{-1}(U')$ and $\bar U:=(\bar U')^\ad=(\psi^{\ad})^{-1} (U)$.
Then, since $p$ is finite, we find
$$
\R_{T'/T} (\psi^{-1}(U'))^\ad = T^\ad \setminus ((p^\ad)^{-1}(T'^\ad \setminus (\psi^\ad)^{-1}(U))).
$$
But now $x\in T^\ad$ is a point, whose $p^\ad$-fibre  $\tilde I$ is contained in $(\psi^\ad)^{-1}(U)$.
This implies $x \in \R_{T'/T} (\psi^{-1}(U'))^\ad$, and hence we finished the proof.
\end{proof}

The previous theorem is a generalization of Bertapelle's result, since we do not require the existence of a formal $R$-model of the rigid space $X'$.

Finally we prove the compatibility of Weil restriction with the functor of adification of locally finite type schemes. (For the notion of this functor cf. \cite{hu3}, 3.8.)
In the following proposition let $F$ denote the forgetful functor from the category of adic spaces into the category of locally ringed spaces.

\begin{proposition}\label{vergleich_adische_algebraische_weilres}
Let $X'$ be a scheme, locally of finite type over $S'$ and let $h \colon S' \longrightarrow S$ be a finite and locally free morphism of schemes.
Furthermore let $T$ be an adic space and let $g \colon F(T) \longrightarrow S$ be a morphism of locally ringed spaces.
We put $T' = T \times_S S'$.
If $X'$ has property (P) for schemes (cf. \cite{blr}, 7.6), the Weil restriction $\R_{T' / T} (X' \times_{S'} T')$ exists in the category of adic spaces over $T$ and coincides with $\R_{S'/S}(X') \times_S T$, the adic space associated to $\R_{S'/S}(X')$.
\end{proposition}
\begin{proof}
We may assume that $S$ and, hence, $S'$ are affine and that $T$ is affinoid.
Since $h$ is finite, so is $h' \colon T' \longrightarrow T$.
We first treat the case where $X'$ is affine.

Let $Z \longrightarrow T$ be a morphism of affinoid adic spaces.
A $T'$-morphism 
$$
Z \times_T T' \longrightarrow X' \times_{S'} T'
$$
is uniquely determined by an $S'$-morphism
$$
F(Z \times_S S') \longrightarrow X'
$$
of locally ringed spaces.
Since the global sections of $Z \times_S S'$ coincide with the ordinary tensor product $\mathcal O_Z(Z) \otimes_{\mathcal{O}_S(S)} \mathcal O_{S'}(S')$ and since $X'$ is affine, such a morphism defines an $S'$-morphism
$$
\Spec \mathcal O_Z(Z) \times_S S' \longrightarrow X'
$$
of schemes in a unique way.
This morphism induces uniquely an $S$-morphism 
$$
\Spec \mathcal O_Z(Z) \longrightarrow \R_{S'/S}(X'),
$$
hence an $S$-morphism $F(Z) \longrightarrow \R_{S'/S}(X')$ of locally ringed spaces, since  $\R_{S'/S}(X')$ is affine.
This again induces by \cite{hu3}, 3.8 a $T$-morphism 
$$
Z \longrightarrow \R_{S'/S}(X') \times_S T
$$
of adic spaces.

If $X'$ is not affine, consider the open covering of $(X_i')_{i\in I}$ of $X'$ by all open affine subschemes $X_i' \subset X'$.
Again let $Z$ be an affinoid $T$-space an let 
$$
a' \colon Z \times_T T' \longrightarrow X' \times_{S'} T'
$$
be a $T'$-morphism of adic spaces.
Then $( X_i' \times_{S'} T' )_{i \in I}$ is an open covering of $X' \times_{S'} T'$ and we put $U_i' := (a')^{-1} ( X_i' \times_{S'} T' )$ and $U_i := \R_{Z \times_T T' / Z} ( U_i' )$ for all $i \in I$.
Using lemma \ref{adische_restriktion_vertraeglich_mit_offenen_immersionen} we obtain that $U_i \subset Z$ is open for all $i \in I$.
Let $z \in Z$ and let $p \colon Z \times_T T' \longrightarrow Z$ be the projection.
Then $(U_i)_{i \in I}$ is an open covering of $Z$ and for all $i \in I$ the restriction  $a_i'$ of $a'$ to $U_i \times_T T'$ factors through $X_i' \times_{S'} T'$.
The affine case yields $T$-morphisms 
$$
U_i \longrightarrow \R_{T'/T} (X_i' \times_{S'} T') = \R_{S'/S} (X_i') \times_{S}T \subset \R_{S'/S} (X') \times_{S} T
$$
corresponding to the restrictions $a_i'$.
Now it is immediate to see that these morphisms glue to a $T$-morphism $a: Z \longrightarrow \R_{S'/S} (X') \times_S T$ which then is the unique $T$-morphism satisfying $\Psi' \circ ( a \times \id_{T'}) = a'$, where $\Psi'$ is the adification of the canonical morphism $\Psi : \R_{S'/S}(X') \times_S S' \longrightarrow X'$.
\end{proof}
	\subsection{Tame base change of adic N\'eron models}

We now want to apply our results on the Weil restriction of adic spaces in order to study the base change of adic N\'eron models.
Let $R$ be a complete discrete valuation ring with fraction field $K$, uniformizer $\pi$ and residue field $k$.

\begin{definition}\label{ad_nm_defi}
Let $X_K$ be a smooth and separated adic space over $\Spa K$.
An \emph{adic N\'eron model} of $X_K$ is a smooth and separated adic space $Y$ over $\Spa R$, together with an open immersion $\tau \colon Y_K \hookrightarrow X_K$, such that
\begin{eqnarray*}
\Hom_{\Spa R} (Z,Y) &\longrightarrow& \Hom_{\Spa K} (Z_K, X_K)\\
a &\longmapsto & \tau \circ (a \times \id_{\Spa K})
\end{eqnarray*}
is bijective for every smooth adic space $Z$ over $\Spa R$.
\end{definition}

For example, the adification $G^\ad$ of the lft N\'eron model $G$ of a split torus $\mathbb G_{m,K}^d$ is an adic N\'eron model of $(\mathbb G_{m,K}^d)^\ad$ (cf. \cite{wahlediss}, 4.1.4).

Let $K'/K$ be a finite Galois extension with group $G$, and let $R'/R$ denote the corresponding extension of the respective discrete valuation rings.
Since $G$ acts on $K'$ via isometries, we obtain a $G$-action on $\Spa R'$, and the morphism $\Spa R' \longrightarrow \Spa R$ is equivariant with the trivial action on $\Spa R$.
Let $X_K$ be a smooth adic space over $\Spa K$ and let $X_{R'}$ be an adic N\'eron model of $X_{K'} = X_K \times_{\Spa K} \Spa K'$.
The group $G$ acts from the right on $X_{K'}$, and the morphism $X_{K'} \longrightarrow \Spa K'$ is $G$-equivariant.
By the universal property of the adic N\'eron model $X_{R'}$, this action extends to an action on $X_{R'}$, where $X_{R'} \longrightarrow \Spa R'$ is equivariant.

We can now define a right action of $G$ on $\R_{\Spa R'/\Spa R} (X_{R'})$, similarly as in  \cite{be}, 2.4.
Namely, for any $g \in G$ and any $T$-valued point $p \in \R_{\Spa R'/\Spa R} (X_{R'})(T)$, we define the $g$-translation to be
\begin{equation}\label{weil_aktion}
p.g := \tau_{X_{R'}}(g) \circ p \circ (\id_T \times \tau_{\Spa R'}(g))^{-1},
\end{equation}
where $\tau_{\Spa R'}(g)$ and $\tau_{X_{R'}}(g)$ denote the $g$-translations on $\Spa R'$ resp. $X_{R'}$.

The functor of fixed points
$$
T \mapsto \Hom_{\Spa R} (T, \R_{\Spa R'/\Spa R}(X_{R'}))^G
$$
for the action of $G$ on $\R_{\Spa R'/\Spa R}(X_{R'})$, where we always want to assume that $X_{R'}$ satisfies property (P), can be represented by an adic space over $\Spa R$.
We now examine the smoothness of this space of fixed points.

Let $A$ be a local noetherian ring with maximal ideal $\mathfrak m$ and let $R \longrightarrow A$ be a local homomorphism of rings.
We denote the $\mathfrak m$-adic completion of $A$ by $\widehat A$.
Let $G$ be a finite group acting from the left on $A$ via continuous translations defined over $R$.
There are two functors of fixed points associated to $A$, namely
\begin{eqnarray}\label{Funktor_F}
F_A \colon \cat{$R$-algebras} & \longrightarrow & \cat{Sets}\\\nonumber
C & \longmapsto & \Hom_R(A,C)^G
\end{eqnarray}
and
\begin{eqnarray}\label{Funktor_F_strich}
F'_A \colon \cat{adic complete $R$-algebras} & \longrightarrow & \cat{Sets}\\\nonumber
C & \longmapsto & \Hom_R^\cont(\widehat{A},C)^G.
\end{eqnarray}
Using the continuity of the action of $G$ on $A$, it is straightforward to show that the following lemma holds:

\begin{lemma}\label{g-aktion-komplettierung}
If $F$ is representable by an $R$-algebra $A^G$ (which is a quotient of $A$), then $F'$ is represented by the $\mathfrak m$-adic completion $\widehat{A^G}$ of $A^G$.
\end{lemma}

\begin{lemma}\label{komplettikompletti}
Let $A$ be an affinoid ring, which is the completion of an affinoid ring $B$.
We denote the canonical morphism by $\tau \colon B \longrightarrow A$. 
Let $\mathfrak p \subset A$ be an open prime ideal and put $\mathfrak q := \tau^{-1}(\mathfrak p)$.
Then the max-adic completions $\widehat{A_\mathfrak p}$ und $\widehat{B_\mathfrak q}$ coincide.
\end{lemma}
\begin{proof}
Let $a$ be an element in $A$ and let $(b_n)_{n \in \mathbb N}$ be a series in $B$ converging towards $a$.
For every $N \in \mathbb N$, there exists an integer $n_0 \ge 0$ satisfying $\tau(b_{n_0}) -a \in \mathfrak p^{N+1}$, hence $\tau(b_{n_0}) + \mathfrak p^{N+1} = a+ \mathfrak{p}^{N+1}$.
It follows that $\tau_N \colon B /\mathfrak q^{N+1} \longrightarrow A / \mathfrak p^{N+1}$ is a bijection.
Hence
$$
B_{\mathfrak q} / \mathfrak{q}^{N+1} B_\mathfrak{q} = (B /\mathfrak{q}^{N+1})_{\overline{\mathfrak{q}}} = (A /\mathfrak{p}^{N+1})_{\overline{\mathfrak{p}}} = A_{\mathfrak p} / \mathfrak{p}^{N+1} A_\mathfrak{p}
$$
for all $N \in \mathbb N$, which implies $\widehat{A_\mathfrak{p}} = \widehat{B_\mathfrak{q}}$.
\end{proof}

\begin{proposition}\label{fixpunkte_glatt}
Let $X$ be an adic space, which is smooth and separated over $\Spa R$ and which has property (P).
Let $G$ be a finite group acting on $X$, such that $X\longrightarrow \Spa R$ is equivariant with the trivial action on $\Spa R$.
Let $n:=\vert G \vert$ be prime to $\cha k$ and assume that the generic fibre $(X^G)_K$ of the space of fixed points is smooth over $\Spa K$.
Then $X^G$ is smooth over $\Spa R$.
\end{proposition}
\begin{proof}
The strategy as follows:
At first we show that the max-adic completions $\widehat{\mathcal O_{X^G,x}}$ in any non-analytic point $x \in X^G$ is formally smooth and, hence, flat over $R$.
Since $X^G$ is locally of finite type over $\Spa R$, it follows that all analytic points of $X^G$ are contained in the generic fibre $\Spa K$.
This yields the $\Spa R$-flatness of $X^G$.
Then we verify that the special fibre of $X^G$ is smooth over $\Spa k$.
Now, using the Jacobi criterion for smoothness of adic spaces (\cite{hu4}, 1.6.9) and the fact that the $\Spa R$-flatness of and affinoid adic space implies the $R$-flatness of its affinoid ring of functions, it is immediate to see that, in our situation, the fibre criterion for smoothness of adic spaces over $\Spa R$ holds.
Hence, our claim is justified.

Now let $x$ be a non-analytic point of $X^G$.
As we assumed $X$ to be smooth over $\Spa R$, \cite{hu4}, (1.7.7), lemma \ref{komplettikompletti} and \cite{egaiv4}, 17.5.3 yield the formal smoothness of $\Spf \widehat{\mathcal O_{X,x}}$ over $\Spf R$.
Since $x$ is $G$-invariant, there is a canonical $G$-action on $\Spf \widehat{\mathcal O_{X,x}}$.
Consider the cartesian square
$$
\xymatrix{
Z \ar[r]^\tau\ar[d]&\Spf \widehat{\mathcal O_{X,x}}\ar[d]\\
\Spf \widehat{\mathcal O_{X,x}} \ar[r]&\Spf R
}
$$
of formal schemes, where $G$ acts on the first factor of $Z$.
Formal smoothness commutes with base change.
Hence $\tau$ is formally smooth.
By \cite{hu1}, 3.6.7 (iii), the ring $\widehat{\mathcal O_{X,x}}$ is noetherian.
We denote the point in $\Spf \widehat{\mathcal O_{X,x}}$ corresponding to the maximal ideal by $s$.
Let $z \in Z$ be the image of $s$ under the diagonal.
Then $k(z) = k(s)$.
The local ring $\widehat{\mathcal O_{Z,z}}$ is formally smooth over $\widehat{\mathcal O_{X,x}}$ and $z$ is fixed under the $G$-action on $Z$.
This yields a $G$-action on $\widehat{\mathcal O_{Z,z}}$.
We now consider the functors $F_{\widehat{\mathcal O_{Z,z}}}$ resp. $F_{\widehat{\mathcal O_{X,x}}}$ (cf. \eqref{Funktor_F}), and we denote the representing $R$-algebras by $\widehat{\mathcal O_{Z,z}}^G$ resp. $\widehat{\mathcal O_{X,x}}^G$.
If we now endow $\widehat{\mathcal O_{Z,z}}^G$ and $\widehat{\mathcal O_{X,x}}^G$ with the quotient topology (which is complete and Hausdorff), it follows that $\widehat{\mathcal O_{Z,z}}^G$ and $\widehat{\mathcal O_{X,x}}^G$ represent the functors $F'_{\widehat{\mathcal O_{Z,z}}}$ resp. $F'_{\widehat{\mathcal O_{X,x}}}$ (cf. \eqref{Funktor_F_strich}).
By definition of the $G$-action on $Z$ we find 
$$
(\widehat{\mathcal O_{Z,z}})^G = (\widehat{\mathcal O_{X,x}})^G \hat \otimes_R \widehat{\mathcal O_{X,x}}.
$$
By faithfully flat descent it now suffices to see that $(\widehat{\mathcal O_{Z,z}})^G$ is formally smooth over $\widehat{ \mathcal O_{X,x}}$ (cf. \cite{wahlediss}, 4.6.4).

As $k(s) = k(z)$, \cite{egaiv1}, 0.19.6.4 yields
$$
\widehat{\mathcal O_{Z,z}} \hat \otimes_{\widehat{\mathcal O_{X.x}}} k(s) = k(s) [[T_1,\ldots,T_d]].
$$
Using \cite{egaiv1}, 0.19.7.1.5 and the flatness of $\widehat{\mathcal O_{X,x}}$ (cf. \cite{egaiv1}, 0.19.7.1), this isomorphism lifts to an isomorphism $B:= \widehat{\mathcal O_{Z,z}} = A [[ T_1,\ldots,T_d]]$, where we have set $A:=\widehat{\mathcal O_{X,x}}$.
By our assumption on $n$ being prime to to $\cha k$, we are in the setting of the proof of \cite{be}, 3.4.
There it is shown that $B^G$ is formally smooth over $A$.
Hence $(\widehat{\mathcal O_{X,x}})^G$ is formally smooth over $R$.
So in order to obtain the formal smoothness of $X^G$ at $x$ it suffices to show that $(\widehat{\mathcal O_{X,x}})^G = \widehat{\mathcal O_{X^G,x}}$ holds.

By lemma \ref{g-aktion-komplettierung},
$$
\Hom_R^\cont ((\widehat{\mathcal O_{X,x}})^G, C) = \Hom_R^\cont(\widehat{\mathcal O_{X,x}},C)^G = \Hom_R^\cont (\widehat{\mathcal O_{X,x}^G},C)
$$
holds for every complete affinoid ring $C$ over $R$, hence $(\widehat{\mathcal O_{X,x}})^G = \widehat{\mathcal O_{X,x}^G}$.
But the $G$-action on $\mathcal O_{X,x}$ is induced by the $G$-action on $X$.
Thus $\mathcal O_{X,x}^G = \mathcal O_{X^G,x}$ and the formal smoothness of $\mathcal O_{X^G,x}$ follows.
In particular, by \cite{egaiv1}, 0.19.7.1, $X^G$ is flat over $\Spa R$.

It remains to show that the special fibre $(X^G)_k$ of $X^G$ is smooth over $\Spa k$.
The space of fixed points $X^G$ is adic over $\Spa R$, hence the special fibre of $X^G$ consists of non-analytic points only.
Let $x \in X^G$.
Since $X$ has property (P), we may assume that $X$ is affinoid, say $X = \Spa A$ for a complete affinoid ring $A$.
Then $A_k$ is a discretely topologized affinoid ring, whose f-adic part $A_k^\tr$ is algebraically smooth over $R$ (cf. \cite{hu4}, 1.6.6 (i)).
By Edixhoven's results, $(\Spec A_k^\tr)^G$ is smooth over $k$.
If we now denote the global sections of $(\Spec A_k^\tr)^G$ by $C$, and define $C^+$ to be the integral closure of $A_k^+$ in $C$, it is easy to check (using again \cite{hu4}, 1.6.6 (i)), that $\Spa C$ represents $(\Spa A_k)^G$ and that $\Spa C$ is smooth over $\Spa k$.
\end{proof}

Let $f \colon X \longrightarrow Y$ be an immersion of adic spaces.
A closed immersion $i \colon \bar X \longrightarrow Y$ is called the \emph{adic closure} of $f$, if $f$ factors through $i$ and if $i$ is minimal with this property.

In the following theorem we compare the adic closure $H$ of $X_K$ in $\R_{R'/R}(X_{R'})$ to the space of fixed points $\R_{R'/R} (X_{R'})$.

\begin{theorem}\label{p_teilt_nicht_grad}
Let $K'/K$ be Galois with group $G$, such that $(\vert G \vert, p) =1$, where $p$ is the residue characteristic.
Let $X_{R'}$ be an adic N\'eron model of $X_{K'}$ with property (P).
Then the adic closure $H$ of $X_K$ in $\R_{R'/R}(X_{R'})$ is the adic N\'eron model of $X_K$.
\end{theorem}
\begin{proof}
At first let us note that the adic closure exists (cf. \cite{wahlediss}, 4.4).
The space of fixed points $\R_{R'/R}(X_{R'})^G$ is a closed subspace of $\R_{R'/R}(X_{R'})$.
The canonical morphism $X_K \longrightarrow \R_{K'/K}(X_{K'})$ is equivariant with the $G$-action on $\R_{K'/K}(X_{K'})$ and the trivial action on $X_K$.
Hence $X_K \subset \R_{R'/R}(X_{R'})^G$.
From the definition of $H$ it follows that $H \subset \R_{R'/R}(X_{R'})^G$ is a closed subspace.

Since the $G$-action of $X_{R'}$ restricts to the canonical $G$-action of $X_{K'}$ and since (in contrast to the formal-rigid setting) Weil restriction commutes with fibres, the equation  
$$
(\R_{R'/R}(X_{R'})^G)_K =\R_{K'/K}(X_{K'})^G
$$
holds.

We now want to verify $\R_{K'/K}(X_{K'})^G = X_K$.
Let $T_K$ be an adic space over $K$.
Consider the canonical morphisms
$$
\Psi \colon \R_{K'/K}(X_{K'}) \times_{\Spa K} \Spa K' \longrightarrow X_{K'}
$$
and
$$
\sigma \colon X_K \longrightarrow \R_{K'/K}(X_{K'}).
$$
Since $\Psi \circ (\sigma \times \id_{\Spa K'}) = \id_{X_{K'}}$, an $\Spa K$-morphism $T_K \longrightarrow X_K$ yields, by definition of the $G$-action on $\R_{K'/K}(X_{K'})$, a $G$-invariant morphism $T_K \longrightarrow \R_{K'/K}(X_{K'})$ of adic spaces over $\Spa K$, by composition with $\sigma$.
There exists a commutative diagram
$$
\xymatrix{
\Hom_K (T_K , X_K) \ar[r] \ar[rd] & \Hom_{K} (T_K,\R_{K'/K} (X_{K'}))^G\ar[d]^\wr\\
&\Hom_{K'} (T_{K'}, X_{K'})^G,
}
$$
where the diagonal map is defined by $a \mapsto a \times \id_{K'}$.
The bijection on the right hand side is induced by $\Psi$.
Hence, according to the definition of the $G$-action on Weil restrictions \eqref{weil_aktion}, it suffices to show that any $\Spa K'$-morphism $a' \colon T_{K'} \longrightarrow X_{K'}$ descends to an $\Spa K$-morphism $a \colon T_K \longrightarrow X_K$, as soon as the diagram
\begin{equation}\label{g-diagramm}
\xymatrix{
T_{K'} \ar[r]^{a'} \ar[d]^{\id \times \tau_g} & X_{K'}\ar[d]^{\id \times \tau_g}\\
T_{K'} \ar[r]^{a'}& X_{K'}
}
\end{equation}
commutes for all $g \in G$.
Here $\tau_g$ denotes the $g$-translation on $\Spa K'$.
We now reduce to the affinoid case.

By assumption, $X_{R'}$ has property (P).
From proposition \ref{eigenschaft_p_steigt_ab} we can conclude that $\R_{K'/K}(X_{K'})$ has property (P), too.
The morphism $\sigma \colon X_K \hookrightarrow \R_{K'/K}(X_{K'})$ is a closed immersion.
Hence $X_K$ has property (P).
We consider the projections $p \colon T_{K'} \longrightarrow T_K$ and $q \colon X_{K'} \longrightarrow X_K$.
For any point $x \in T_K$ it follows from the finiteness of $p$ that there exists an open affinoid neighborhood $U_K \subset X_K$ of $q(a'(p^{-1}(x)))$.
Then $U_{K'}:=U_K \times_{\Spa K} \Spa K'$ is an open affinoid neighborhood of $a'(p^{-1}(x))$ and
$$
\tilde V_K := T_K \setminus p( T_{K'} \setminus (a')^{-1}(U_{K'}))
$$
is an open neighborhood of $x$ in $T_K$.
We choose an open affinoid neighborhood $V_K \subset \tilde V_K$ of $x$.
Now the open affinoid $V_K\times_{\Spa K} \Spa K'$ is contained in $(a')^{-1}(U_{K'})$.

Hence we may assume $T_K = \Spa B$ and $X_K = \Spa A$, for complete affinoid rings $A$ and $B$.
We write
\begin{eqnarray*}
&p_1^* \colon K' \longrightarrow K' \otimes_K K', \quad \alpha' \longmapsto \alpha' \otimes 1\\
\text{ and } &p_2^* \colon K' \longrightarrow K' \otimes_K K',\quad  \alpha' \longmapsto 1 \otimes \alpha',
\end{eqnarray*}
and we denote the field extension $K'/K$ by $p^* \colon K \longrightarrow K'$ and put $q^* := p_1^* \circ p^* = p_2^* \circ p^*$.
We claim that the diagram
$$
\xymatrix{
\Hom_{K}^\cont (A,B) \ar[r] &\Hom_{K'}^\cont (p^*A, p^*B) \ar@<2pt>[r]^{p_1} \ar@<-2pt>[r]_{p_2}&\Hom_{K' \otimes_K K'}^\cont(q^*A, q^*B)
}
$$
of sets of continuous morphisms of affinoid rings is exact.
This follows immediately from Galois descent of $K$-algebras.
Namely, the $K'$-morphism $(a')^* \colon A \otimes_K K' \longrightarrow B \otimes_K K'$ of affinoid rings descends to a $K$-homomorphism $a^* \colon A \longrightarrow B$ of ordinary $K$-algebras.
It suffices to explain that this homomorphism is in fact a continuous morphism of affinoid rings.
But the fact that $a^*$ is compatible with rings of integral elements follows immediately from the construction of the affinoid tensor product.
The injective map $B \longrightarrow B \otimes_K K'$ is strict (\cite{hu1}, 2.3.33 (ii), resp. \cite{hu1}, 2.2.15 (ii)).
Hence $a^*$ is continuous.
This means that
$$
(\R_{R'/R}(X_{R'})^G)_K = \R_{K'/K}(X_{K'})^G = X_K.
$$

In particular, $H_K = (\R_{R'/R}(X_{R'})^G)_K$, and we claim $H = \R_{R'/R}(X_{R'})^G$.

By proposition \ref{fixpunkte_glatt}, $\R_{R'/R}(X_{R'})^G$ is $R$-flat and $H \subset \R_{R'/R} (X_{R'})^G$ is a closed subspace.
Let $(Y_j)_{j \in J}$ be an open affinoid covering of $\R_{R'/R}(X_{R'})^G$, and let $j \in J$.
Then $\mathcal O_{Y_j}(Y_j)$ is flat over $R$.
The intersection $Y_j \cap H$ is a closed subspace of $Y_j$.
Now the $R$-flatness of $\mathcal O_{Y_j} (Y_j)$ allows us to conclude  $Y_j \cap H = Y_j$  from $(Y_j \cap H)_K = (Y_j)_K$ for all $j \in J$, hence $H = \R_{R'/R}(X_{R'})^G$.
Finally, proposition \ref{fixpunkte_glatt} yields the smoothness of $\R_{R'/R}(X_{R'})^G$.
\end{proof}

In analogy to \cite{bert}, section 2, we may now deduce the following result from what we have proven so far.

\begin{corollary}\label{tame_base_change}
Let $X_K$ be a smooth adic space over $\Spa K$ with adic N\'eron model $X_R$ over $\Spa R$ and let $K'/K$ be a tamely ramified field extension, such that $X_{K'}:=X_K \times_{\Spa K} \Spa K'$ admits an adic N\'eron model $X_{R'}$ over $\Spa R'$ with property (P).
Then $X_R$ is isomorphic to the adic closure of $X_K$ in $\R_{R'/R}(X_{R'})$.
\end{corollary}

	\section{Weil restriction of Berkovich spaces}
	\subsection{The affinoid case}\label{berko_affinoide}

In this section let $K$ be a field, which is complete with respect to a non-archimedean, non-trivial valuation $\vert \cdot \vert$.

The following definition is a generalization of the notion of the spectral value of a monic polynomial with coefficients in a (strictly) $K$-affinoid algebra.

\begin{definition}
Let $A$ be a $K$-affinoid Banach algebra and consider a polynomial
$$
p(z) = z^n + c_1 z^{n-1} + \ldots + c_n \in A [z].
$$
If $\rho$ denotes the spectral radius on $A$, $\sigma (p) := \max_{i=1,\ldots, n} \rho( c_i )^{1/i}$ is called the \emph{spectral value of $p$}.
\end{definition}

\begin{lemma}\label{spektralwertgleichspektralradius}
Let $\varphi \colon A \longrightarrow B$ be a bounded homomorphism of $K$-affinoid Banach algebras and let $B$ be a finite free $A$-module with $A$-basis $e_1,\ldots,e_n \in B$.
Write $b = \sum_{i=1}^n b_i e_i \in B$, $b_1,\dots,b_n \in B$ and consider the characteristic polynomial $\chi_b(z) \in A [z]$ of $b$ over $A$.
Then $$\rho(b) = \sigma(\chi_b).$$

In particular, if $\chi_b (z) = z^n + c_1 z^{n-1} + \ldots + c_n$, with coefficients $c_1,\dots,c_n \in A$ and if $r$ is a positive real number, then $\rho(b) \le r$ iff $\rho( c_i ) \le r^i$ for all $i=1,\ldots,n$.
\end{lemma}
\begin{proof}
We reduce to the strictly $K$-affinoid case, where the statement is true by lemma 1.6 in \cite{bert}.
Let $r=(r_1,\ldots,r_n)$ be a system of positive real numbers, such that $A \hat \otimes_K K_r$ and, hence, $B \hat \otimes_K K_r$, is strictly $K_r$-affinoid.
The vertical arrows in the commutative diagram
\[
 \begin{xy}
  \xymatrix{ A \hat \otimes_K K_r \ar[r]&B \hat \otimes_K K_r\\
A \ar[u] \ar[r]&B, \ar[u]
}
 \end{xy}
\]
are isometries.
For any $b \in B$, 
\[
\rho(b) = \lim_{n \to \infty} \sqrt[n]{ \Vert b^n \Vert }= \lim_{n \to \infty} \sqrt[n]{ \Vert (b \hat \otimes 1)^n \Vert} = \rho_r (b \hat \otimes 1),
\]
where $\rho_r$ denotes the spectral radius on $B \hat \otimes_K K_r$.
The characteristic polynomial $\chi_{b \hat \otimes 1}$ of $b \hat \otimes 1$ over $A \hat \otimes_K K_r$ coincides with the polynomial obtained by applying $A \longrightarrow A \hat \otimes_K K_r$ to the coefficients of $\chi_b$.
Lemma 1.6 in \cite{bert} shows $\sigma_r (\chi_{b \hat \otimes 1}) = \rho_r (b \hat \otimes 1)$, and we conclude $\sigma_r (\chi_{b \hat \otimes 1}) = \sigma (\chi_b)$.
\end{proof}

Let $\Phi_K$ be a system of $K$-affinoid Berkovich spaces and let $h \colon M(B) \longrightarrow M(A)$ be a finite, free morphism of rank $n$ of $\Phi_K$-affinoid Berkovich spaces. (A morphism $h$ of $\Phi_K$-affinoid Berkovich spaces will be called \emph{free}, if the associated bounded homomorphism $h^*$ of $K$-Banach algebras is a free ring homomorphism.)
Let $e_1,\ldots,e_n \in B$ be an $A$-basis of $B$.
We now construct the Weil restriction of $\mathbb{B}^1_{B,r}$ over $B$ with radius $r > 0, r \in \mathbb{R}$ with respect to the base change $h$.

Let $X= (x_1,\ldots,x_n)$ and $Y = (y_1,\ldots,y_n)$ be finite systems of indeterminates.
The characteristic polynomial of $\sum_{i=1}^n x_i e_i \in B[X]$ will be denoted by $\chi (z) \in A[X] [z]$.
We write
$$
\chi(z) = z^n + c_1(X) z^{n-1} + \ldots + c_n(X) 
$$
with coefficients $c_1(X),\ldots,c_n(X) \in A[X]$.
For all $s > 0, s \in \mathbb{R}$ we define
\begin{equation}
 C_{s,r} := A \langle s^{-1} X, r^{-1}y_1, \ldots, r^{-n}y_n \rangle / (y_1 - c_1(X),\ldots,y_n - c_n(X)), \label{C_s_r-Definition}
\end{equation}
where we have used the shorthand notation ``$s^{-1} X$'' instead of ``$s^{-1}x_1,\ldots,s^{-1}x_n$''.
If $\pi$ denotes the canonical projection from $A \langle s^{-1} X, r^{-1}y_1, \ldots, r^{-n}y_n \rangle$ to $C_{s,r}$, we obtain
\[
\rho ( \pi(c_i(X)))= \rho (\pi(y_i)) \le \rho (y_i) = r^{i}
\]
for all $i=1,\ldots,n$.

If $\Phi_K$ coincides with the system of all strictly $K$-affinoid Banach algebras and if $r \in \vert K^* \vert$ (say $r = \vert \alpha \vert$), the Berkovich space $\mathbb B^1_{B,r}$ is associated to the rigid disc $\Sp B \langle r^{-1} x \rangle$ over $\Sp B$ with radius $r$.
In this case, $C_{s,r}$, $s=\vert \pi^\lambda \vert$, $\lambda \in \mathbb N$ is also strictly $K$-affinoid and the $K$-algebras $C_{s,r}$, $C_{\lambda,\{\alpha\}}$ and the corresponding algebras defined by Bertapelle in \cite{bert} are isomorphic.
This follows from the following simple fact from linear algebra:

If $b \in B$ and $c_i(b) \in A$, $1 \le i \le n$ the $i$-th coefficient of the characteristic polynomial $\chi_b(z)$ of $b$ over $A$, then
\[
 c_i(ab) = a^i c_i(b) \qquad \text{for all }i=1,\ldots,n, \text{ and all } a\in A,
\]
since, for the characteristic polynomial $\chi_{ab}(z)$ of $ab$, we find $\chi_{ab}(z) = a^n \chi_b (z/a)$.

Now we glue the $\Phi_K$-affinoid spaces $M(C_{s,r})$.
We write 
$$
M_0 (C_{s,r}) := \{ z \in M(C_{s,r}) ; \vert x_i \vert_z < s \text{ for all } i=1,\ldots,n\}.
$$

\begin{lemma}
Let $r \in \mathbb R$, $r>0$.
For each $s >0$, $M_0 (C_{s,r})$ is an open subspace of $M(C_{s,r})$ and there are open immersions 
$$
M_0 (C_{s,r}) \hookrightarrow M_0(C_{s',r})
$$
for all $s' \ge s$.
\end{lemma}
\begin{proof}
Let $s' \ge s$.
Then $\rho (x_i) \le s \le s'$ for all $i=1,\ldots,n$.
By \cite{berk1}, 2.1.5 there is a unique bounded homomorphism $\iota_{s',s}^* \colon C_{s',r} \longrightarrow C_{s,r}$ satisfying $x_i \mapsto x_i$ and $y_i \mapsto y_i$ for all $i=1,\ldots,n$.
From $M(C_{s,r}) = \{ z \in M (C_{s',r}) ; \vert x_i \vert_z \le s \text{ for all } i=1,\ldots,n\}$ it follows that $M(C_{s,r}) \subset M(C_{s',r})$ is an affinoid subdomain.
By definition of the topology on $M(C_{s,r})$, $M_0 (C_{s,r}) \subset M (C_{s,r})$ is open for all $s>0$, and $M_0(C_{s,r})$ is isomorphic to an open subspace of $M_0(C_{s',r})$.
\end{proof}

Now we can glue the quasi-affinoid spaces $M_0 (C_{s,r})$, $s >0$.
The result is a $\Phi_K$-analytic space over $M(A)$, which we will denote by $\varinjlim_{s>0} M_0 (C_{s,r})$.

\begin{proposition}\label{berko_kreisscheibe_res}
With the above notations, there exists a functorial isomorphism
$$\R_{M(B) / M(A)} (\mathbb{B}^1_{B,r}) \cong \varinjlim_{s>0} M_0 (C_{s,r}).$$
\end{proposition}
\begin{proof}
We only need to verify the universal property for $\Phi_K$-affinoid spaces $M(D)$ over $M(A)$.
For $s>0$ we define
\[
\Psi_s^* \colon B \langle r^{-1}x \rangle \longrightarrow C_{s,r} \otimes_A B
\]
by putting $\Psi_s^*(x) := \sum_{j=1}^n x_je_j$.
By lemma \ref{spektralwertgleichspektralradius}, this is a well-defined, bounded homomorphism of $K$-Banach algebras.
We set $\Psi_s := M(\Psi_s^*)$ and $\Psi := \varinjlim_{s>0} \Psi_s$.
Since $M(D)$ is compact, one can easily verify that the map
\[
\Hom_{M(A)} (M(D),\varinjlim_{s>0} M_0 ( C_{s,r})) \longrightarrow \Hom_{M(B)} (M(D) \times_{M(A)} M(B), \mathbb{B}^1_{B,r}),
\]
given by $a \mapsto \Psi \circ (a \times \id_{M(B)})$, is a bijection.
\end{proof}

From the compatibility of Weil restriction with products we can deduce that the Weil restriction of arbitrary polydiscs over $B$ exists.
We formulate this result in the following proposition.

\begin{proposition}\label{berko_polykreis_darstellbar}
Let $X=(x_1,\ldots,x_d)$ be a finite system of indeterminates and $r=(r_1,\ldots,r_d)$ a finite system of positive real numbers.
Then the Weil restriction of
\[
\mathbb{B}^d_{B,r} := M(B \langle r^{-1}X \rangle ) =M(B \langle r_1^{-1}x_1, \ldots, r_d^{-1}x_d \rangle )
\]
exists, and
\begin{eqnarray*}
\R_{M(B)/M(A)} (\mathbb{B}^d_{B,r}) &=& \prod_{i=1}^d \R_{M(B)/M(A)} (\mathbb{B}^1_{B,r_i})
\\
&=&\varinjlim_{s>0} M_0 \left( E_{s,r} \right),
\end{eqnarray*}
where we have written $E_{s,r}:=\widehat{\bigotimes}_{i=1,\ldots,d} C_{s,r_i}$.
Here, $M_0 (E_{s,r})$ is defined as follows:
For all $i=1,\ldots,d$ let $X_i = (x_{i1},\ldots, x_{in})$ be the system of the first $n$ coordinates on $C_{s,r_i}$ (cf. \eqref{C_s_r-Definition} on page \pageref{C_s_r-Definition}).
We denote the open subspace of $M (E_{s,r})$, on which the values of the coordinates $x_{ij}$, $i=1,\ldots,d$, $j=1,\ldots ,n$ are strictly smaller than $s$, by $M_0 (E_{s,r})$.
\end{proposition}

Now we prove the compatibility of Weil restriction with closed immersions.

\begin{proposition}\label{berko_restriktion_affinoid}
Let $D$ be a $\Phi_K$-affinoid Banach algebra over $B$.
Then the Weil restriction $\R_{M(B)/M(A)} (M(D))$ is a good $\Phi_K$-analytic Berkovich space over $M(A)$.
\end{proposition}
\begin{proof}
Let $\alpha \colon B \langle r^{-1} X \rangle \longrightarrow D$ be an admissible epimorphism of $K$-affinoid Banach algebras, where $X = (x_1,\ldots,x_d)$ is a finite system of variables and $r = (r_1,\ldots,r_d)$ a finite system of real numbers.
We set $\mathfrak a := \ker \alpha$.
Then $\mathfrak a$ is finitely generated by \cite{berk1}, 2.1.3, say $\mathfrak a=(f_1,\ldots,f_k)$ for certain elements $f_1,\ldots,f_k \in B \langle r^{-1} X \rangle$.
Let $e_1, \ldots, e_n$ be an $A$-basis of $B$ and put $E_{s,r} := \widehat{\bigotimes}_{i=1,\ldots,d} C_{s,r_i}$.
For each $s >0$ we have the canonical maps 
\[
\Psi_s^* \colon B \langle r^{-1}X \rangle \longrightarrow E_{s,r} \otimes_A B.
\]
For each $j=1,\ldots,k$ we write 
$$
\Psi_s^* (f_j) = \sum_{i=1}^n f_{ijs} e_i,\quad f_{1js},\ldots,f_{njs} \in E_{s,r}
$$
and denote by $\mathfrak a_s^\co$ the ideal generated by $f_{ijs}$, $i=1,\ldots,n$, $j=1,\ldots,k$ in $E_{s,r} \otimes_A B$.
Keeping the notations from proposition \ref{berko_polykreis_darstellbar}, we put $M_0 (E_{s,r} / \mathfrak a^\co_s) := M_0(E_{s,r})\cap M(E_{s,r}/\mathfrak a^\co_s)$, hence
$$
M_0 (E_{s,r} / \mathfrak a^\co_s) = \{ z \in M(E_{s,r} / \mathfrak a^\co_s) ; \vert \overline{x_{i1}} \vert_z,\ldots,\vert \overline{x_{in}} \vert_z  < s \text{ for all } i=1,\ldots,d \},
$$
where $\overline{x_{ij}} \in E_{s,r} / \mathfrak a^\co_s$ denotes the image of the $j$-th coordinate function $x_{ij}$ of the $i$-th factor $C_{s,r_i}$ of $E_{s,r}$.
Similar to the adic setting, one easily verifies that there exists a bijection
\[
\Hom_{M(A)} (T, \varinjlim_{s>0} M_0 (E_{s,r} / \mathfrak a^\co_s)) \longrightarrow \Hom_{M(B)} (T \times_{M(A)} M(B), M(D))
\]
for all $\Phi_K$-affinoid spaces $T$ over $M(A)$.
\end{proof}
	\subsection{Topological properties of finite morphisms of Berkovich spaces}
In this section we derive a few basic topological properties of finite morphisms of $\Phi_K$-analytic spaces, where  $\Phi_K$ will always denote a system of $K$-affinoid spaces as defined in \cite{berk2}.
The proofs of lemma \ref{endlich_folgt_hausdorffsch} and proposition \ref{endliche_K_analytisch_folgt_kompakt} are based on a correspondence with V. Berkovich.

\begin{lemma}\label{endlich_folgt_hausdorffsch}
A finite morphism $f \colon X \longrightarrow Y$ of $\Phi_K$-analytic Berkovich spaces is Hausdorff.
\end{lemma}
\begin{proof}
Let $x_1, x_2 \in X$ satisfying $f (x_1) = f (x_2) =:y$ and $x_1 \neq x_2$.
Then there are $\Phi$-affinoids $V_1, \ldots, V_n$ in $Y$ such that $V_1 \cup \ldots \cup V_n$ is a neighborhood of $y$ in $Y$.
The sets $U_i := f^{-1}(V_i)$, $1 \le i \le n$ have the property that $x_1, x_2 \in \bigcup_{i=1}^n U_i$ holds and that $\bigcup_{i=1}^n U_i$ is a neighborhood of $x_1$ and $x_2$.
Since $f$ is finite, the subsets $U_i$, $1 \le i \le n$ are $K$-affinoid and, in particular, Hausdorff.
Then for each $i=1,\ldots,n$ there are open neighborhoods $V_{i1}$ of $x_1$ and $V_{i2}$ of $x_2$, such that $V_{i1}\cap V_{i2} \cap U_i$ is empty.
Thus $\bigcap_{i=1}^n V_{ij}$, $j=1,2$ are disjoint, open neighborhoods of the $x_j$.
\end{proof}

\begin{proposition}\label{endliche_K_analytisch_folgt_kompakt}
A finite morphism $f \colon X \longrightarrow Y$ of $K$-analytic Berkovich spaces is a compact map of topological spaces.
\end{proposition}
\begin{proof}
Let $A \subset Y$ be compact, $B = f^{-1}(A)$.
Then, by lemma \ref{endlich_folgt_hausdorffsch}, $B$ is Hausdorff.
Let $(V_i)_{i \in I}$ be an open covering of $B$.
A set $U \subset Y$ will be called a \emph{special neighborhood of $y \in Y$}, if $U$ is contained in an open, Hausdorff neighborhood of $y$, and if $U$ is a neighborhood of $y$, which can be written as $U = U_1 \cup \ldots \cup U_n$, where $U_1,\ldots,U_n$ are $K$-affinoid domains, such that $y \in U_1 \cap \ldots \cap U_n$.
Special neighborhoods are compact, and every $y \in Y$ has a fundamental system of special neighborhoods.
It now suffices to show:
\\

\emph{$(\ast)$ Every $y \in A$ has a special neighborhood $U'$, such that $f^{-1}(U')$ is contained in the union of a finite subcollection $(V_j)_{j \in J}$, $J \subset I$ of $(V_i)_{i \in I}$.}
\\

Namely, we cover $A$ by special neighborhoods $U_y'$, $y \in A$, satisfying property $(\ast)$.
A finite number of $U_y'$'s will do, and we write $A \subset \bigcup_{i=1}^r U_i'$ with special neighborhoods $U_1',\ldots, U_r'$.
It follows that there exists finite subsets $J_1 ,\ldots , J_r \subset I$ such that $f^{-1} (U_i') \subset \bigcup_{j \in J_i} V_j$ for all $i=1,\ldots,r$.
Hence $(V_j)_{j \in J_i, i=1,\ldots,r}$ is a finite subcovering of $(V_i)_{i \in I}$ which still covers $B$.

In order to verify $(\ast)$ we first observe that, due to the finiteness of $f$, for all $y \in A$ there is a neighborhood $V$ of the fibre $f^{-1}(y)$ of type $V = V_{i_1} \cup \ldots \cup V_{i_m}$, $i_1,\ldots,i_m \in I$, $m \ge 1$.
Let $U$ be a special neighborhood of $y$, say $U = U_1 \cup \ldots \cup U_n$, $y \in U_1 \cap \ldots \cap U_n$, $n \ge 1$.
Then $f^{-1}(U)$ is Hausdorff by lemma \ref{endlich_folgt_hausdorffsch}.
This implies that the restrictions $f\vert_{f^{-1}(U)} \colon f^{-1}(U) \longrightarrow U$ of $f$ to $f^{-1}(U)$ are compact.
Namely, if $C \subset U$ is a compact subset of the Hausdorff set $U$, then $C$ is closed in $U$, and the finiteness of $f$ implies that $f^{-1}(C)$ can be written as a finite union of compact spaces in a Hausdorff space, which is again compact.
The set $V \cap f^{-1}(U)$ is an open neighborhood of $f^{-1}(y)$ in $f^{-1}(U)$.
But now $f\vert_{f^{-1}(U)} \colon f^{-1} (U) \longrightarrow U$ is a compact map and $U$ is compact.
Hence $f\vert_{f^{-1}(U)}$ is a closed map (cf. \cite{bourt}, §10.3, Proposition 7).
Then $\tilde U := U \setminus f( f^{-1}(U) \setminus V \cap f^{-1}(U))$ is an open neighborhood of $y$, whose preimage is contained in $V \cap f^{-1} (U)$.
Since the special neighborhoods of $y$ constitute a fundamental system of neighborhoods of $y$, we can choose a special neighborhood $U' \subset \tilde U$ of $y$, such that $f^{-1} (U') \subset V \cap f^{-1}(U) \subset V = V_{i_1} \cup \ldots \cup V_{i_m}$.
\end{proof}

Using the canonical functor $\BERK \longrightarrow \KAN$ and \cite{berk2}, lemma 1.3.7, we can conclude:

\begin{corollary}
A finite morphism $f \colon X \longrightarrow Y$ of $\Phi_K$-analytic spaces is compact.
\end{corollary}

\begin{corollary}\label{berko_topologisch_abgeschlossen}
Let $f \colon X \longrightarrow Y$ be a finite morphism of $\Phi_K$-analytic Berkovich spaces.
Then $f$ is a closed map of topological spaces.
\end{corollary}
\begin{proof}
By \cite{berk2}, 1.2.4 (iii), every point $y \in Y$ possesses a compact neighborhood $U_y \subset Y$.
As we have seen, $f^{-1}(U_y)$ is compact for all $y \in Y$.
Using \cite{bourt}, I, §10.1 Proposition 3 b), it suffices to show that the restrictions
$$
f_y \colon f^{-1}(U_y) \longrightarrow U_y,\qquad y \in Y,
$$
of $f$ are universally closed, and this follows from \cite{bourt}, I, §10.2 Corollaire 2.
\end{proof}	
	\subsection{A criterion for good spaces}\label{berko_gute_abschnitt}
There are two ways to glue Berkovich spaces from local ones.
The first is to glue along open subspaces, the second is to glue along analytic subdomains, which are closed but not necessarily open.
In the second case one has to require that, locally, only finitely many overlaps occur.
In the first case this limitation does not apply.
However, if we glue Weil restrictions of quasi-affinoid subspaces of our global object, we cannot expect a criterion for the representability of the Weil restriction for the non-good case.
Hence, we treat the good and the non-good case separately.

Let $h \colon S' \longrightarrow S$ be a finite and locally free morphism of $\Phi_K$-analytic Berkovich spaces and let $U'\subset S'$ be an open $\Phi$-analytic subdomain.
Then $U:=S \setminus h (S' \setminus U')$ is open in $S$.
If we know that $U$ is a $\Phi$-analytic subdomain, then it is clear that $U$ represents $\R_{S'/S}(U')$.
For example, this is the case if we assume $\Phi$ to be dense.

\begin{lemma}\label{quasiaffinoidexistiert}
Let $\Phi_K$ be a dense system of $K$-affinoid spaces.
Let $X'$ be a $\Phi_K$-analytic Berkovich space over $S'$ and let $h \colon S' \longrightarrow S$ be finite and locally free.
If $\R_{S'/S}(X')$ can be represented by a $\Phi_K$-analytic Berkovich space over $S$, then, for all open subspaces $U' \subset X'$, $\R_{S'/S}(U')$ is an open subspace of $\R_{S'/S}(X')$.
In particular, the Weil restriction of quasi-affinoid $K$-spaces exists.
\end{lemma}
\begin{proof}
We set $T:= \R_{S'/S}(X')$ and $T' := T \times_S S'$.
Then $U' \times_{X'} T'$ is isomorphic to an open subdomain $W'$ of $T'$ and 
$$
\R_{T'/T}(W') = \R_{S'/S}(U') \times_{\R_{S'/S}(X')} T = \R_{S'/S} (U'),
$$
by the compatibility of Weil restriction with base change.
\end{proof}

\begin{theorem}\label{berko_darstellbar_gut}
Let $\Phi$ be a dense system of affinoid $K$-spaces, $(X',\mathcal A',\tau')$ a $\Phi_K$-analytic Berkovich space and let $h \colon S' \longrightarrow S$ be a finite free morphism of $\Phi_K$-affinoid spaces.
Assume that for all finite subsets $I \subset X'$ there exists a $\Phi$-affinoid subdomain $U' \subset X'$, which is a neighborhood of $I$ in $X'$.
Then $\R_{S'/S}(X')$ is a good $\Phi_K$-analytic Berkovich space over $S$.
\end{theorem}
\begin{proof}
By assumption, $X'$ is good.
Hence, $(U_V)_{V \in \widehat \tau'}$ is a covering of $X'$, where we denoted the interior of $V \in \widehat \tau'$ by $U_V$.
By lemma \ref{quasiaffinoidexistiert}, all functors $\R_{S'/S} (U_V)$, $V \in \widehat \tau'$ are representable, and we obtain a $\Phi_K$-analytic space $\R$ by glueing the representing $\Phi_K$-analytic spaces along the canonical glueing data.
Also, the canonical $S$-morphisms $\Psi_V \colon \R_{S'/S} (U_V) \times_S S' \longrightarrow U_V$, $V \in \widehat \tau'$ glue to an $S$-morphism $\Psi \colon \R \times_S S' \longrightarrow X'$.

We eventually prove that $\Psi$ defines a bijection
\begin{eqnarray*}
\Theta \colon \Hom_S (Z,\R) &\longrightarrow& \Hom_{S'} (Z \times_S S',X')\\
a &\longmapsto & \Psi \circ ( a \times \id_{S'})
\end{eqnarray*}
for all $K$-affinoid spaces $Z$ over $S$.

Let $a' \colon Z \times_S S' \longrightarrow X'$ be an $S'$-morphism.
Then $p \colon Z \times_S S' \longrightarrow Z$ is finite.
Hence $p^{-1} (z)$ is finite, and there exists a $K$-affinoid neighborhood $V_z$ of $a'(p^{-1}(z))$ in $X'$.
In particular $a' ( p^{-1}(z))$ is contained in $U_{V_z}$, and  $W_z := \R_{Z \times_S S'/Z} ((a')^{-1}(U_{V_z}))$ is an open neighborhood of $z$ in $Z$.
Given the restriction $a'_z \colon W_z \times_S S' \longrightarrow U_{V_z}$ of $a'$, there is a unique $S$-morphism $a_z \colon W_z \longrightarrow \R_{S'/S} (U_{V_z})$ satifying $a'_z = \Psi_{V_z} \circ ( a_z \times \id_{S'})$.
The family $(W_z)_{z \in Z}$ yields a covering $Z$, along with morphisms $a_z \colon W_z \longrightarrow \R_{S'/S} (U_{V_z})$.
It is immediate to see that the morphisms $a_x$ agree on overlaps, and, hence, yield an $S$-morphism $a \colon W \longrightarrow Z$, which satisfies $\Theta (a) = a'$.

In the same way the injectivity of $\Theta$ can be proven.
Let $z \in Z$.
Due to property (P) there exists a $V \in \widehat \tau'$ such that
$$
\Theta (a_1)(p^{-1}(z)) = \Theta (a_2)(p^{-1}(z)) \subset U_V
$$
and we set $W:= \R_{Z \times_S S'/Z} (\Theta(a_1)^{-1}(U_V))$.
Then $W$ is an open neighborhood of $z$ in $Z$.
By definition of $W$, 
$$
\Theta(a_1) \vert_{W \times_S S'} =  \Theta(a_2) \vert_{W \times_S S'}
$$
factors through $U_V \subset X'$.
We claim that $a_1 \vert_W$ and $a_2 \vert_W$ factor through $\R_{S'/S} (U_V)$.
Let $w \in W$ and let $L$ be the completion of the residue field in $w$.
We denote the canonical morphism $M(L) \longrightarrow Z \hat \otimes_K L$ again by $w$.
Then by \cite{berk2}, 1.4, $p^{-1}(w)$ is homeomorphic to the $\Phi_L$-analytic space $M(L) \times_{S_L} S'_L$, where $S_L := S \hat \otimes_K L$ and $S'_L := S' \hat \otimes_K L$.
Let $V' \in \widehat \tau'$ such that $a_1(w) \in \R_{S'/S} (U_{V'})$.
From the representability of $\R_{S'_L/S_L}(U_{V'} \hat \otimes_K L)$ it follows that
$$
a_{1,w} := a_1 \circ w \colon M(L) \longrightarrow \R_{S'/S}(U_{V'}) \hat \otimes_K L = \R_{S'_L/S_L}(U_{V'} \hat \otimes_K L)
$$
induces an $S'_L$-morphism
$$
a_{1,w}' \colon M(L) \times_{S_L} S'_L \longrightarrow U_{V'} \hat \otimes_K L
$$
of $\Phi_L$-analytic Berkovich spaces.
By definition of $\Psi$, $a_{1,w}'$ and $(\Theta(a_1) \hat \otimes \id_L) \circ (w \times \id_{S'_L})$ coincide, and $a_{1,w}'$ factors through $(U_{V'} \cap U_V) \hat \otimes_K L$.
Again, from the representability of $\R_{S'_L/S_L}(U_{V'} \hat \otimes_K L)$ it follows that $a_{1,w}$ already factors through $\R_{S'_L/S_L}((U_{V'} \cap U_V) \hat \otimes_K L) \subset \R_{S'/S}(U_V) \hat \otimes_K L$.
Hence, the image of $a_1 \vert_W$ is contained in $\R_{S'/S}(U_V)$.
Moreover, $a_2 (W) \subset \R_{S'/S}(U_V)$.
But then, similarly as in the proof of theorem \ref{adisch_darstellbar}, we can conclude $a_1=a_2$.
\end{proof}
	\subsection{A criterion for not necessarily good space}

In order to deal with non-good spaces, we now need to glue along closed $\Phi$-analytic domains.
For this it is again necessary to determine the Weil restriction of a closed $\Phi$-analytic domain $V$, which is already contained in the base $S'$.
Since $h$ is an open map, $W:=S \setminus h(S' \setminus V)$ is a closed subset of $S$.
Once we can show that this subset has the structure of a $\Phi_K$-analytic space, it is clear that $W$ represents the Weil restriction of $V$ with respect to $h$.
In the following we only deal with the case where $\Phi_K$ consists of all $K$-affinoid spaces.

\begin{lemma}\label{berk_vertraeglich_analytisch}
Let $h \colon S' \longrightarrow S$ be a finite and locally free morphism of $K$-analytic Berkovich spaces and let $V \subset S'$ be a closed $K$-analytic subdomain of $S'$.
Then $W := S \setminus h(S' \setminus V)$ is a closed $K$-analytic subdomain of $S$.
\end{lemma}
\begin{proof}
It suffices to prove that $W$ is analytic.
For this we may assume that both $S$ and $S'$ are affinoid and where $h$ is finite and \emph{free}, i.e. on global sections we have a finite and free morphism of rings.
Furthermore, we may assume that $V$ is a special domain in $S'$.

We now reduce to the case where $S$, $S'$ and $V$ are strictly $K$-analytic.
There exists a finite system $r$ of positive real numbers, such that $S \hat \otimes_K K_r$ (and, hence, $S' \hat \otimes_K K_r$) and $V \hat \otimes_K K_r$ are strictly $K_r$-analytic.
Consider the cartesian square
$$
\xymatrix{
S' \hat \otimes_K K_r \ar[r]^{h'} \ar[d]^{\pi'} & S \hat \otimes_K K_r \ar[d]^\pi\\
S' \ar[r]^h & S.
}
$$
Then $\pi^{-1}(W)$ consists of all those points in $S \hat \otimes_K K_r$ whose $h'$-fibers are contained in $V \hat \otimes_K K_r$.

Assume our claim is true in the strict case.
Then
$$
W':=S \hat \otimes_K K_r \setminus h' ( S' \hat \otimes_K K_r \setminus V \hat \otimes_K K_r)
$$
is an analytic closed subset of $S \hat \otimes_K K_r$.
Then, using  \cite{duc}, 2.4, $W'$ can be written as a finite union of rational domains.
From the proof of \cite{duc}, 2.4 it follows that then also $W$ can be written as such a finite union of rational domains.
Hence $W \subset S$ is analytic.

It remains to verify the claim under the assumption that $S$, $S'$ and $V$ are strictly $K$-affinoid.
Consider the commutative diagram
$$
\xymatrix{
M(B) \ar[r]^h &M(A)\\
\Sp B \ar[u] \ar[r]^{h_0}&\Sp A. \ar[u]
}
$$
The finiteness of $h$ yields $W \cap \Sp A = \widetilde W$, where we have set
$$
\widetilde W := \Sp A \setminus h_0(\Sp B \setminus (V \cap \Sp B)).
$$
Since $V \cap \Sp B \subset \Sp B$ is admissible open and quasi-compact, $\widetilde W \subset \Sp A$ is also admissible open and quasi-compact by \cite{bert}, 1.12. 
Consider the functor $X \mapsto X_0$, sending a strictly $K$-analytic Hausdorff space $X$ to its set $X_0$ of rigid points (which can be canonically endowed with the structure of a rigid $K$-space, cf. \cite{berk2}, 1.6).
By \cite{berk2}, 1.6.1, the functor $X \mapsto X_0$ is fully faithful and induces an equivalence of categories between the category of paracompact strictly $K$-analytic Berkovich spaces and the category of quasi-separated rigid $K$-spaces which admit an admissible covering of finite type.
In particular, associated to $\widetilde W$ there is a $K$-analytic Berkovich space $X$ such that $X_0 = \widetilde W$.
For all $K$-affinoid Berkovich spaces $Z$ we have equations
\begin{eqnarray*}
\Hom_{M(A)} (Z,X) &=& \Hom_{\Sp A}(Z_0, \widetilde W)\\
&=&\Hom_{\Sp B}(Z_0 \times_{\Sp A} \Sp B, V_0)\\
&=&\Hom_{M(B)}(Z \times_{M(A)} M(B),V)\\
&=&\R_{M(B)/M(A)}(V)(Z).
\end{eqnarray*}
By construction of $X$ (cf. proof of theorem 1.6.1 in \cite{berk2}) and by lemma 1.6.2 (i) in \cite{berk2} we see that $X \subset M(A)$ is a closed analytic subspace.
First we want to show that $X$ represents $\R_{M(B)/M(A)}(V)$ for the category of all $K$-analytic Berkovich spaces.
After that it only remains to check $W = X$ in order to finish the proof.

Let $Z$ be a $K$-analytic space over $S$.
We have seen that $\R_{S'/S}(V)$ is representable by $X$ for the category of strictly $K$-analytic spaces over $S$.
Hence, if $Z$ is strictly $K$-analytic, a morphism $Z \times_S S' \longrightarrow S'$ factors through $V \subset S'$ if and only if $Z \longrightarrow S$ is contained in $X$.

If $Z$ is not strictly $K$-analytic, there exists a $K$-affinoid field $K_r$, such that $Z \hat \otimes_{K} K_r$ is strictly $K_r$-affinoid.
A morphism $Z \longrightarrow S$ factors through $X \subset S$ if and only if $Z \hat \otimes_K K_r \longrightarrow S \hat \otimes_K K_r$ factors through $X \hat \otimes_K K_r$.
The latter is equivalent to the fact that $(Z \times_S S') \hat \otimes_K K_r \longrightarrow S' \hat \otimes_K K_r$ factors through $V \hat \otimes_K K_r$, since $X \hat \otimes_K {K_r}$ represents $\R_{S' \hat \otimes_K K_r/S \hat \otimes_K K_r} (V \hat \otimes_K K_r)$ for the category of strictly $K_r$-analytic spaces.
Finally, this last statement is equivalent to the statement that the image of $Z \times_S S' \longrightarrow S'$ is contained in $V$.
Thus $X$ represents $\R_{S'/S}(V)$ for the category of all $K$-analytic spaces over $S$.
If we regard a point $x \in S$ as an $M(L)$-valued point of $S \hat \otimes_K L$, where $L$ is again the completion of the residue field at $x$, the representability of $\R_{S'/S}(V)$ by $X$ yields that $x$ is contained in $X$ iff $h^{-1}(x)$ is contained in $V$, hence $W=X$.
\end{proof}

By replacing ``open $K$-analytic domain'' by ``closed $K$-analytic domain'' in the proof of lemma \ref{quasiaffinoidexistiert} and by using lemma \ref{berk_vertraeglich_analytisch} instead of corollary \ref{berko_topologisch_abgeschlossen}, we can conclude:

\begin{proposition}\label{berk_vertraeglich_analytisch_satz}
Let $h \colon S' \longrightarrow S$ be a finite and locally free morphism of $K$-analytic Berkovich spaces.
Let $(X',\mathcal A',\tau')$ be a $K$-analytic Berkovich space over $S'$ and let $U' \subset X'$ be a closed $K$-analytic domain in $X'$.
Assume $\R_{S'/S} (X')$ is representable by a $K$-analytic Berkovich space over $S$.
Then $\R_{S'/S}(U')$ is representable by a closed $K$-analytic domain of $\R_{S'/S}(X')$.
\end{proposition}

\begin{theorem}\label{berko_schlechte_res_darstellbar}
Let $h \colon S' \longrightarrow S$ be a finite free morphism of $K$-affinoid Berkovich spaces.
Let $(X',\mathcal A', \tau')$ be a $K$-analytic Berkovich space over $S'$ such that $\tau'$ is locally finite.
Assume that for all finite subsets $I \subset X'$ there exists a $V \in \tau'$ such that $I \subset V$.
Then $\R_{S'/S}(X')$ exists.
\end{theorem}
\begin{proof}
By assumption and by the compactness of $V$, $\{ W \in \tau' ; W \cap V \neq \emptyset \}$ is finite for each $V \in \tau'$.
For every $V \in \tau'$ the Weil restriction $\R_{S'/S}(V)$ exists, cf. proposition \ref{berko_restriktion_affinoid}.
If $W \in \tau'$, then $W \cap V \subset V$ is a closed $K$-analytic domain.
By proposition \ref{berk_vertraeglich_analytisch_satz} the Weil restriction $\R_{S'/S}(W \cap V)$ is a closed $K$-analytic domain of $\R_{S'/S}(V)$.
The set $\{ W \in \tau' ; \R_{S'/S} (W \cap V) \neq \emptyset \}$ is finite.
Then the $K$-analytic spaces $\R_{S'/S}(V)$, $V \in \tau'$ can be glued to a $K$-analytic Berkovich space $\R$, which then has a covering by closed $K$-analytic domains $R_V := \R_{S'/S}(V)$, $V \in \tau'$.

We now consider  the net of all domains $W \subset \R$, which admit a $V \in \tau'$ such that $W \subset R_V$ and which correspond under the isomorphism $\R_{S'/S}(V) \longrightarrow R_V$ to an affinoid domain in $\R_{S'/S}(V)$.
In particular, each point $x \in \R$ has a neighborhood $W_1 \cup \ldots \cup W_n$, $x \in W_1 \cap \ldots \cap W_n$ such that $W_i  \subset R_{V_i}$ for all $i=1,\dots,n$.
It follows that $R_{V_1} \cup \ldots \cup R_{V_n}$ is a neighborhood of $x$ in $\R$.
Hence $(R_V)_{V \in \tau'}$ is a quasi-net on $\R$.
By \cite{berk2}, 1.3.2, the canonical morphisms $\Psi_V \colon \R_{S'/S}(V) \times_S S' \longrightarrow V \subset X'$, $V \in \tau'$ induce a morphism $\Psi \colon \R  \times_S S' \longrightarrow X'$.
Now let $Z$ be a $K$-affinoid Berkovich space over $S$.
Similarly to the proof of theorem \ref{berko_darstellbar_gut} we now show that
\begin{eqnarray*}
\Theta \colon \Hom_{S} ( Z, \R) &\longrightarrow & \Hom_{S'} (Z \times_S S', X')\\
a & \longmapsto &\Psi \circ (a \times \id_{S'})
\end{eqnarray*}
is a bijection.

Let $a' \colon Z' \longrightarrow X'$ be an $S'$-morphism, let $x \in Z$ and let $p \colon Z' \longrightarrow Z$ denote the projection onto the first factor, where $Z' := Z \times_S S'$.
Then there exists a $V_x \in \tau'$ containing $a'(p^{-1}(x))$ such that $p^{-1}(x) \subset (a')^{-1}(V_x)$, this means that $x$ is an element of $W_x:=\R_{Z'/Z}((a')^{-1}(V_x))$.
By lemma \ref{berk_vertraeglich_analytisch} we know that $(W_x)_{x \in Z}$ is a covering of $Z$ by closed $K$-analytic domains in $Z$.
Since $Z'$ is compact, $((a')^{-1}(V_x))_{x \in Z}$ and, hence, $(W_x)_{x \in Z}$ admits a finite subcovering.
Thus $(W_x)_{x \in Z}$ is a quasi-net on $Z$.
Restricting $a'$ to $W_x \times_S S' \subset Z'$, we obtain $S$-morphisms 
$$
a_x \colon W_x \longrightarrow \R_{S'/S}(V_x) \hookrightarrow \R.
$$
Using proposition \ref{berk_vertraeglich_analytisch_satz}, the $a_x$ agree on the intersections
$$
W_x \cap W_y = \R_{Z'/Z}((a')^{-1}(V_x \cap V_y))
$$
for all $x,y \in Z$ and thus, according to \cite{berk2}, 1.3.2, give rise to an $S$-morphism $a \colon Z \longrightarrow \R$ satisfying $\Psi \circ (a \times \id_{S'}) = a'$.

Now consider $S$-morphisms $a_1, a_2 \colon Z \longrightarrow \R$ such that 
$$a' = \Psi \circ (a_1 \times \id_{S'}) = \Psi \circ (a_2 \times \id_{S'}).$$
Let $x \in Z$.
Then there exists an affinoid $V_x \in \tau'$ such that $a' (p^{-1}(x)) \subset V_x$.
It follows that $W_x := \R_{Z'/Z}((a')^{-1}(V_x))$ is a closed $K$-analytic domain in $Z$, which contains $x$.
Since both $a_1 \vert_{W_x}$ and $a_2 \vert_{W_x}$ factor through $\R_{S'/S} (V_x)$, we obtain, in the same manner as in the proof of theorem \ref{berko_darstellbar_gut}, that $a_1 \vert_{W_x} = a_2 \vert_{W_x}$ holds, which concludes the proof.
\end{proof}	
	\subsection{Weil restriction of non-quasi-projective schemes}\label{non_qp_berk}

Consider a finite field extension $K'/K$, where still $K$ is a field which is complete with respect to a non-trivial,  non-archimedean valuation.
If $X'$ is a quasi-projective $K'$-scheme, the Weil restriction $\R_{K'/K}(X')$ of $X'$ exists as a $K$-scheme (cf. \cite{blr}, 7.6/4).

We can use our results on the Weil restriction of good Berkovich spaces in order to see that the Weil restriction of $X'$ with respect to $K'/K$ already exists as a good Berkovich space over $K$, if $X'$ is separated and locally of finite type over $K'$.

\begin{proposition}
Let $K'/K$ be a finite field extension, let $X'$ be a separated $K'$-scheme, locally of finite type over $K'$, and let $(X')^\an$ be the $K'$-analytic Berkovich space associated to $X'$.
Then $\R_{K'/K}((X')^\an)$ is a good Berkovich space over $K$.
\end{proposition}
\begin{proof}
Of course, we use theorem \ref{berko_darstellbar_gut}.
Since $X'$ is separated, $(X')^\an$ is a good Hausdorff Berkovich space over $K'$.
Let $I \subset X'$ be a finite set.
Then for each $x \in I$ there exists an open neighborhood $U_x$ such that the $U_x$, $x \in I$ are pairwise disjoint sets.
As $(X')^\an$ is good, so is $U_x$.
Hence, for each $x \in I$ there exists an affinoid neighborhood $V_x$ contained in $U_x$.
But then $V:=\bigcup_{x \in I} V_x$ is an affinoid neighborhood of $I$ in $X$.
\end{proof}

We conjecture that $\R_{K'/K}((X')^\an)$ is the Berkovich space associated to an algebraic space over $K$ (in the sense of \cite{ct}).			
	
	\bibliography{main} 

\providecommand{\bysame}{\leavevmode\hbox to3em{\hrulefill}\thinspace}
\providecommand{\MR}{\relax\ifhmode\unskip\space\fi MR }
\providecommand{\MRhref}[2]{%
  \href{http://www.ams.org/mathscinet-getitem?mr=#1}{#2}
}
\providecommand{\href}[2]{#2}
\begin{thebibliography}{10}

\bibitem{berk1}
Vladimir~G. Berkovich, \emph{Spectral theory and analytic geometry over
  non-{A}rchimedean fields}, Mathematical Surveys and Monographs, vol.~33,
  American Mathematical Society, Providence, RI, 1990. \MR{MR1070709
  (91k:32038)}

\bibitem{berk2}
\bysame, \emph{\'{E}tale cohomology for non-{A}rchimedean analytic spaces},
  Inst. Hautes \'Etudes Sci. Publ. Math. \textbf{78} (1993), 5--161 (1994).
  \MR{MR1259429 (95c:14017)}

\bibitem{bert}
Alessandra Bertapelle, \emph{Formal {N}\'eron models and {W}eil restriction},
  Math. Ann. \textbf{316} (2000), no.~3, 437--463. \MR{MR1752779 (2001g:14038)}

\bibitem{bgr}
Siegfried Bosch, Ulrich G{\"u}ntzer, and Reinhold Remmert,
  \emph{Non-{A}rchimedean analysis}, Grundlehren der Mathematischen
  Wissenschaften, vol. 261, Springer-Verlag, Berlin, 1984. \MR{MR746961
  (86b:32031)}

\bibitem{blr}
Siegfried Bosch, Werner L{\"u}tkebohmert, and Michel Raynaud, \emph{N\'eron
  models}, Ergebnisse der Mathematik und ihrer Grenzgebiete (3), vol.~21,
  Springer-Verlag, Berlin, 1990. \MR{MR1045822 (91i:14034)}

\bibitem{bourt}
Nicolas Bourbaki, \emph{Topologie g\'en\'erale}, Diffusion C.C.L.S., Paris,
  1971.

\bibitem{ct}
Brian Conrad and Michael Temkin, \emph{Non-archimedean analytification of
  algebraic spaces}, Preprint, \verb+http://arxiv.org/abs/0706.3441v1+, 2007.

\bibitem{duc}
Antoine Ducros, \emph{Parties semi-alg\'ebriques d'une vari\'et\'e alg\'ebrique
  {$p$}-adique}, Manuscripta Math. \textbf{111} (2003), no.~4, 513--528.
  \MR{MR2002825 (2004m:14122)}

\bibitem{be}
Bas Edixhoven, \emph{N\'eron models and tame ramification}, Compositio Math.
  \textbf{81} (1992), no.~3, 291--306. \MR{MR1149171 (93a:14041)}

\bibitem{egaiv1}
Alexander Grothendieck, \emph{\'{E}l\'ements de g\'eom\'etrie alg\'ebrique.
  {IV}. \'{E}tude locale des sch\'emas et des morphismes de sch\'emas. {I}},
  Inst. Hautes \'Etudes Sci. Publ. Math. \textbf{20} (1964), 259. \MR{MR0173675
  (30 \#3885)}

\bibitem{egaiv4}
\bysame, \emph{\'{E}l\'ements de g\'eom\'etrie alg\'ebrique. {IV}. \'{E}tude
  locale des sch\'emas et des morphismes de sch\'emas. {IV}}, Inst. Hautes
  \'Etudes Sci. Publ. Math. \textbf{32} (1967), 361. \MR{MR0238860 (39 \#220)}

\bibitem{hu1}
Roland Huber, \emph{Bewertungsspektrum und rigide {G}eometrie}, Regensburger
  Mathematische Schriften, vol.~23, Universit\"at Regensburg Fachbereich
  Mathematik, Regensburg, 1993. \MR{MR1255978 (95c:32036)}

\bibitem{hu2}
\bysame, \emph{Continuous valuations}, Math. Z. \textbf{212} (1993), no.~3,
  455--477. \MR{MR1207303 (94e:13041)}

\bibitem{hu3}
\bysame, \emph{A generalization of formal schemes and rigid analytic
  varieties}, Math. Z. \textbf{217} (1994), no.~4, 513--551. \MR{MR1306024
  (95k:14001)}

\bibitem{hu4}
\bysame, \emph{\'{E}tale cohomology of rigid analytic varieties and adic
  spaces}, Aspects of Mathematics, E30, Friedr. Vieweg \& Sohn, Braunschweig,
  1996. \MR{MR1734903 (2001c:14046)}

\bibitem{wahlediss}
Christian Wahle, \emph{Die {T}echnik der {W}eil-{R}estriktion f\"{u}r {H}uber-
  und {B}erkovich-{R}\"{a}ume und {A}nwendungen zu adischen
  {N}\'{e}ron-{M}odellen}, Ph.D. thesis, Westf\"{a}lische
  Wilhelms-Universit\"{a}t, 2009,
  \verb+http://wwwmath.uni-muenster.de/u/cwahle/docs/diss_wahle.pdf+.

\bibitem{weil}
Andr{\'e} Weil, \emph{Adeles and algebraic groups}, Progress in Mathematics,
  vol.~23, Birkh\"auser Boston, Mass., 1982. \MR{MR670072 (83m:10032)}

\end{thebibliography}
\end{document}